\documentclass[12pt]{amsart}
\usepackage{amssymb}
\usepackage{amsmath}
\usepackage{amsfonts}

\newtheorem{theorem}{Theorem}[section]
\theoremstyle{plain}

\newtheorem{corollary}[theorem]{Corollary}

\newtheorem{definition}[theorem]{Definition}
\newtheorem{example}[theorem]{Example}

\newtheorem{lemma}[theorem]{Lemma}

\newtheorem{proposition}[theorem]{Proposition}
\newtheorem{remark}[theorem]{Remark}

\newtheorem{Observation}[theorem]{Observation}
\numberwithin{equation}{section}

\input{tcilatex}

\begin{document}
\title[Finite Group Crossed-Products]{Irreducible Representations of
C*-crossed products by Finite Groups}
\author{Alvaro Arias}
\address{Department of Mathematics\\
University of Denver\\
Denver CO 80208}
\email{aarias@math.du.edu}
\urladdr{http://www.math.du.edu/\symbol{126}aarias}
\thanks{Part of this research took place when the first author visited Texas
A\&M University to participate in the Workshop in Analysis and Probability.
He expresses his appreciation for their hospitality.}
\author{Fr\'{e}d\'{e}ric Latr\'{e}moli\`{e}re}
\address{Department of Mathematics\\
University of Denver\\
Denver CO 80208}
\email{frederic@math.du.edu}
\urladdr{http://www.math.du.edu/\symbol{126}frederic}
\date{February 10, 2010}

\begin{abstract}
We describe the structure of the irreducible representations of crossed
products of unital C*-algebras by actions of finite groups in terms of
irreducible representations of the C*-algebras on which the groups act. We
then apply this description to derive a characterization of irreducible
representations of crossed-products by finite cyclic groups in terms of
representations of the C*-algebra and its fixed point subalgebra. These
results are applied to crossed-products by the permutation group on three
elements and illustrated by various examples.
\end{abstract}

\maketitle

\section{Introduction}

What is the structure of irreducible representations of C*-crossed-products $%
A\rtimes _{\alpha }G$ of an action $\alpha $ of a finite group $G$ on a
unital C*-algebra $A$? Actions by finite groups provide interesting
examples, such as quantum spheres \cite{Bratteli91,Bratteli92} and actions
on the free group C*-algebras \cite{CLat06}, among many examples, and have
interesting general properties, as those found for instance in \cite%
{Rieffel80}. Thus, understanding the irreducible representations of their
crossed-products is a natural inquiry, which we undertake in this paper.

\bigskip After we wrote this paper, we are shown \cite{Takesaki67} where Takesaki provides in  a detailed description of the irreducible representations
of $A\rtimes G$ when $G$ is a locally group acting on a type I C*-algebra $A$ and the action is assumed
{\it smooth}, as defined in \cite[Section 6]{Takesaki67}. Our paper takes a different road, though with some important intersections we were not aware of originally.
Both our paper and \cite{Takesaki67} make use of the Mackey machinery and the structure of the commutant of the image of irreducible representations of the crossed-products. However, since our C*-algebras are not assumed to be type I, and in general the restriction of an irreducible representation of $A\rtimes G$ to $A$ does not
lead to an irreducible representation of $A$, we need a different approach than \cite{Takesaki67}. The main tool we use for this purpose is the impressive result proven in \cite{Hoegh-Krohn81} that for ergodic actions of compact groups on unital C*-algebras, spectral subspaces are finite dimensional. As a consequence, we can analyze irreducible representations of finite group crossed-products on arbitrary unital C*-algebras with no condition on the action of the group on the spectrum of the C*-algebra. More formally, we restrict the assumption on the group and relax it completely on the C*-algebra and the action compared to \cite[Theorem 7.2]{Takesaki67}.

\bigskip Our research on this topic was initiated in a paper of Choi and the
second author \cite{Latremoliere06} in the case where $G=\mathbb{Z}_{2},$
i.e. for the action of an order two automorphism $\sigma $ on a C*-algebra $%
A $. In this situation, all irreducible representations of $A\rtimes
_{\sigma }\mathbb{Z}_{2}$ are either \emph{minimal}, in the sense that their
restriction to $A$ is already irreducible, or are regular, i.e. induced by a
single irreducible representation $\pi $ of $A$ such that $\pi $ and $\pi
\circ \sigma $ are not equivalent. In this paper, we shall answer the
question raised at the beginning of this introduction for any finite group $%
G $. Thus, we suppose given any action $\alpha $ of $G$ on a unital
C*-algebra $A$. In this general situation, we show that for any irreducible
representation $\Pi $ of $A\rtimes _{\alpha }G$ on some Hilbert space $%
\mathcal{H}$, the group $G$ acts ergodically on the commutant $\Pi
(A)^{\prime }$ of $\Pi (A)$, and thus, by a theorem of Hoegh-Krohn, Landstad
and Stormer \cite{Hoegh-Krohn81}, we prove that $\Pi (A)^{\prime }$ is
finite dimensional. We can thus deduce that there is a subgroup $H$ of $G$
such that $\Pi $ is constructed from an irreducible representation $\Psi $
of $A\rtimes _{\alpha }H$, with the additional property that the restriction
of $\Psi $ to $A$ is the direct sum of finitely many representations all
equivalent to an irreducible representation $\pi $ of $A$. In addition, the
group $H$ is exactly the group of elements $h$ in $G$ such that $\pi $ and $%
\pi \circ \alpha _{h}$ are equivalent. The canonical unitaries of $A\rtimes
_{\alpha }G$ are mapped by $\Pi $ to generalized permutation operators for
some decomposition of $\mathcal{H}$. This main result is the matter of the
third section of this paper.

\bigskip When $G$ is a finite cyclic group, then we show that the
representation $\Psi $ is in fact minimal and obtain a full characterization
of irreducible representations of $A\rtimes _{\alpha }G$. This result can
not be extended to more generic finite groups, as we illustrate with some
examples. In addition, the fixed point C*-subalgebra of $A$ for $\alpha $
plays a very interesting role in the description of minimal representations
when $G$ is cyclic. We investigate the finite cyclic case in the fourth
section of this paper.

\bigskip We then apply our work to the case where $G$ is the permutation
group $\mathfrak{S}_{3}$ on three elements $\left\{ 1,2,3\right\} $. It is
possible again to fully describe all irreducible representations of any
crossed-product $A\rtimes _{\alpha }\mathfrak{S}_{3}$, and we illustrate all
the cases we can encounter by examples. This matter is discussed in the last
section of this paper.

\bigskip We start our paper with a section on generalities on
crossed-products of C*-algebras by finite groups, including a result on a
characterization of irreducible regular representations. This section also
allows us to set some of our notations. We now fix some other notations
which we will use recurrently in this paper. Given a Hilbert space $\mathcal{%
H}$ which we decompose as a direct sum $\mathcal{H}=\mathcal{H}_{1}\oplus
\ldots \oplus \mathcal{H}_{m}$ of Hilbert subspaces, we shall write an
operator $T$ on $\mathcal{H}$ as an $m\times m$ matrix whose $(i,j)$-entry
is the operator $p_{i}Tp_{j}$ where $p_{1},\ldots ,p_{m}$ are the orthogonal
projections from $\mathcal{H}$ onto respectively $\mathcal{H}_{1},\ldots ,%
\mathcal{H}_{m}$. If $t_{1},\ldots ,t_{m}$ are operators on, respectively, $%
\mathcal{H}_{1},\ldots ,\mathcal{H}_{m}$, then the diagonal operator with
entries $t_{1},\ldots ,t_{m}$ will be denoted by $t_{1}\oplus \ldots \oplus
t_{m}$, i.e. $\oplus _{j=1}^{m}t_{j}\left( \xi _{1},\ldots ,\xi _{m}\right)
=\left( t_{1}\xi _{1},\ldots ,t_{m}\xi _{m}\right) $ for all $\left( \xi
_{1},\ldots ,\xi _{i}\right) \in \mathcal{H}_{1}\oplus \ldots \oplus 
\mathcal{H}_{m}$. If $\pi _{1},\ldots ,\pi _{m}$ are representations of some
C*-algebra $A$ acting respectively on $\mathcal{H}_{1},\ldots ,\mathcal{H}%
_{m}$, then the representation $\pi _{1}\oplus \ldots \oplus \pi _{m}$ of $A$
on $\mathcal{H}$ is defined by $\left( \pi _{1}\oplus \ldots \oplus \pi
_{m}\right) (a)=\pi _{1}(a)\oplus \ldots \oplus \pi _{m}(a)$ for all $a\in A$%
. The identity operator of $\mathcal{H}$ will be denoted by $1_{\mathcal{H}}$
or simply $1$ when no confusion may occur. More generally, when an operator $%
t$ on a Hilbert space $\mathcal{H}$ is a scalar multiple $\lambda 1_{%
\mathcal{H}}$ ($\lambda \in \mathbb{C}$) of the identity of $\mathcal{H}$ we
shall simply denote it by $\lambda $ and omit the symbol $1_{\mathcal{H}}$
when appropriate.

We shall denote by $f_{|E_{0}}$ the restriction of any function $%
f:E\longrightarrow F$ to a subset $E_{0}$ of $E$. The set $\mathbb{T}$ is
the unitary group of $\mathbb{C}$, i.e. the set of complex numbers of
modulus 1.

\section{Crossed-Product By Finite Groups}

\bigskip In this paper, we let $A$ be a unital C*-algebra and $\alpha $ an
action on $A$ of a finite group $G$ by *-automorphisms. A covariant
representation of $\left( A,\alpha ,G\right) $ on a unital C*-algebra $B$ is
a pair $\left( \pi ,V\right) $ where $\pi $ is a *-homomorphism from $A$
into $B$ and $V$ is a group homomorphism from $G$ into the unitary group of $%
B$ such that for all $g\in G$ and $a\in A$ we have $V(g)\pi (a)V\left(
g^{-1}\right) =\pi \circ \alpha _{g}(a)$. The crossed-product C*-algebra $%
A\rtimes _{\alpha }G$ is the universal C*-algebra among all the C*-algebras
generated by some covariant representation of $\left( A,\alpha ,G\right) $.
In particular, $A\rtimes _{\alpha }G$ is generated by a copy of $A$ and
unitaries $U^{g}$ for $g\in G$ such that $U^{gh}=U^{g}U^{h}$, $%
U^{g^{-1}}=\left( U^{g}\right) ^{\ast }$ and $U^{g}aU^{g^{-1}}=\alpha
_{g}(a) $ for all $g,h\in G$ and $a\in A$. The construction of $A\rtimes
_{\alpha }G$ can be found in \cite{Pedersen79} and is due originally to \cite%
{Zeller-Meier68}.

\bigskip By universality, crossed-products by finite groups have a very
simple form which we now describe.

\begin{proposition}
\label{Presentation}Let $G$ be a finite group of order $n$ and write $%
G=\left\{ g_{0},\ldots ,g_{n-1}\right\} $ with $g_{0}$ the neutral element
of $G$. Let $\sigma $ be the embedding of $G$ in the permutation group of $%
\left\{ 0,\ldots ,n-1\right\} $ given by $\sigma _{g}(i)=j$ if and only if $%
gg_{i}=g_{j}$ for all $i,j\in \left\{ 0,\ldots ,n-1\right\} $ and $g\in G$.
We now define $V_{g}$ to be the matrix in $M_{n}(A)$ whose $(i,j)$ entry is
given by $1_{A}$ if $\sigma _{g}(i)=j$ and $0$ otherwise, i.e. the tensor
product of the permutation matrix for $\sigma _{g}$ and $1_{M_{n}}(A)$. Let $%
\psi :A\longrightarrow M_{n}(A)$ be the *-monomorphism:%
\begin{equation*}
\psi :a\in A\longmapsto \left[ 
\begin{array}{cccc}
a &  &  &  \\ 
& \alpha _{g_{1}}(a) &  &  \\ 
&  & \ddots &  \\ 
&  &  & \alpha _{g_{n-1}}(a)%
\end{array}%
\right] \text{.}
\end{equation*}%
Then $A\rtimes _{\alpha }G$ is *-isomorphic to $\oplus _{g\in G}\psi
(A)V_{g} $. In particular: 
\begin{equation*}
\dbigoplus_{g\in G}AU^{g}=A\rtimes _{\alpha }G\text{.}
\end{equation*}
\end{proposition}

\begin{proof}
The embedding of $G$ into permutations of $G$ is of course the standard
Cayley Theorem. We simply fix our notations more precisely so as to properly
define our embedding $\psi $. A change of indexing of $G$ simply correspond
to a permutation of the elements in the diagonal of $\psi $ and we shall
work modulo this observation in this proof. For $b\in M_{n}(A)$ we denote by 
$b_{i,i^{\prime }}$ its $(i,i^{\prime })$-entry for $i,i^{\prime }\in
\left\{ 1,\ldots ,n\right\} $.

An easy computation shows that:%
\begin{equation*}
V_{g}\psi (a)V_{g^{-1}}=\psi \left( \alpha _{g}(a)\right)
\end{equation*}%
and $V_{g}V_{h}=V_{gh}$ for all $g,h\in G$ and $a\in A$. Therefore, by
universality of $A\rtimes _{\alpha }G$, there exists a (unique)
*-epimorphism $\eta :A\rtimes _{\alpha }G\twoheadrightarrow \oplus _{g\in
G}\psi (A)V_{g}$ such that $\eta _{|A}=\psi $ and $\eta (U^{g})=V_{g}$ for $%
g\in G$. Our goal is to prove that $\eta $ is a *-isomorphism.

First, we show that $\oplus _{g\in G}AU^{g}$ is closed in $A\rtimes _{\alpha
}G$.

Let $\left( a_{m}^{0},\ldots ,a_{m}^{n-1}\right) _{m\in \mathbb{N}}$ in $%
A^{n}$ such that $\left( \sum_{j=0}^{n-1}a_{m}^{j}U^{g_{j}}\right) _{m\in 
\mathbb{N}}$ is a convergent sequence in $A\rtimes _{\alpha }G$. Now:%
\begin{equation*}
\eta \left( \sum_{j=0}^{n-1}a_{m}^{j}U^{g_{j}}\right) =\sum_{j=0}^{n-1}\psi
(a_{m}^{j})V_{g}\text{.}
\end{equation*}%
By definition, we have $\sigma _{g_{j}}(0)=i$ for all $i\in \left\{ 0,\ldots
,n-1\right\} $. Let $j\in \left\{ 0,\ldots ,n-1\right\} $. Then $%
V_{j+1,1}^{g_{j}}=1_{A}$ and $V_{j+1,1}^{g_{i}}=0$ for all $i\in \left\{
0,\ldots ,n-1\right\} \backslash \left\{ j\right\} $. Hence, $\left( \eta
\left( \sum_{j=1}^{n}a_{m}^{j}U^{g_{j}}\right) \right) _{1,j+1}=a_{m}^{j}$
for all $m\in \mathbb{N}$. Since $\eta $ is continuous, and so is the
canonical projection $b\in M_{n}(A)\longrightarrow b_{1,j+1}\in A$, we
conclude that $\left( a_{m}^{j}\right) _{m\in \mathbb{N}}$ converges in $A$.
Let $a^{j}\in A$ be its limit. Then $\left( a_{m}^{0},\ldots
,a_{m}^{n-1}\right) _{m\in \mathbb{N}}$ converges in $A^{n}$ to $\left(
a^{0},\ldots ,a^{n-1}\right) $. Thus, $\left(
\sum_{j=0}^{n-1}a_{m}^{j}U^{g_{j}}\right) _{m\in \mathbb{N}}$ converges to $%
\sum_{j=0}^{n-1}a^{j}U^{g_{j}}\in \oplus _{g\in G}AU^{g}$ and thus $\oplus
_{g\in G}AU^{g}$ is closed in $A\rtimes _{\alpha }G$. Since $\oplus _{g\in
G}AU^{g}$ is dense in $A\rtimes _{\alpha }G$ by construction, we conclude
that $A\rtimes _{\alpha }G=\oplus _{g\in G}AU^{g}$.

Now, we show that $\eta $ is injective. Let $c\in A\rtimes _{\alpha }G$ such
that $\eta (c)=0$. Then there exists $a_{0},\ldots ,a_{n-1}\in A$ such that $%
c=\sum_{j=0}^{n-1}a_{j}U^{g_{j}}$. Let $j\in \left\{ 0,\ldots ,n-1\right\} $%
. Then $\eta (c)=0$ implies that $\eta (c)_{j+1,1}=a_{j}=0$ for all $j\in
\left\{ 0,\ldots ,n-1\right\} $ and thus $c=0$. So $\eta $ is a
*-isomorphism and our proof is concluded.
\end{proof}

\bigskip As we will focus our attention on the crossed-products by finite
cyclic groups in the fourth section of this paper and Proposition (\ref%
{Presentation}) is particularly explicit in this case, we include the
following corollary:

\begin{corollary}
Let $\sigma $ be an automorphism of order $n$ of a unital C*-algebra $A$.
Then $A\rtimes _{\sigma }\mathbb{Z}_{n}$ is *-isomorphic to:%
\begin{equation*}
\left\{ 
\begin{array}{c}
\left[ 
\begin{array}{ccccc}
a_{1} & a_{2} & a_{3} & \cdots & a_{n} \\ 
\sigma (a_{n}) & \sigma (a_{1}) & \sigma (a_{2}) & \sigma (a_{3}) &  \\ 
\sigma ^{2}(a_{n-1)} & \sigma ^{2}(a_{n}) & \sigma ^{2}(a_{1}) & \ddots & 
\ddots \\ 
\vdots & \ddots & \ddots & \ddots & \sigma ^{n-2}(a_{2}) \\ 
\sigma ^{n-1}(a_{2}) & \sigma ^{n-1}(a_{3}) & \cdots & \sigma ^{n-1}(a_{n})
& \sigma ^{n-1}(a_{1})%
\end{array}%
\right] \in M_{n}(A) \\ 
a_{1},\ldots ,a_{n}\in A%
\end{array}%
\right\}
\end{equation*}%
where $U^{1}$ mapped to $\left[ 
\begin{array}{cccc}
0 & 1 &  & 0 \\ 
\vdots & 0 & \ddots &  \\ 
0 & \vdots & \ddots & 1 \\ 
1 & 0 & \cdots & 0%
\end{array}%
\right] $ and $A$ is embedded diagonally as $a\in A\mapsto \left[ 
\begin{array}{cccc}
a &  &  &  \\ 
& \sigma (a) &  &  \\ 
&  & \ddots &  \\ 
&  &  & \sigma ^{n-1}(a)%
\end{array}%
\right] $. In particular, $A\rtimes _{\sigma }\mathbb{Z}_{n}=A\oplus
AU^{1}\oplus \ldots \oplus AU^{n-1}$.
\end{corollary}

\begin{proof}
Simply write $\mathbb{Z}_{n}=\left\{ 0,\ldots ,n-1\right\} $ so that: 
\begin{equation*}
V_{1}= \left[ 
\begin{array}{cccc}
0 & 1 &  & 0 \\ 
\vdots & 0 & \ddots &  \\ 
0 & \vdots & \ddots & 1 \\ 
1 & 0 & \cdots & 0%
\end{array}
\right] \text{.}
\end{equation*}
The result is a direct computation of $\oplus _{k=0}^{n-1}\psi (A)\left(
V_{1}\right) ^{k}$.
\end{proof}

\bigskip We now turn our attention to the irreducible representations of $%
A\rtimes _{\alpha }G$. Proposition (\ref{Presentation}) suggests that we
construct some representations from one representation of $A$ and the left
regular representation of $G$. Of particular interest is to decide when such
representations are irreducible. We will use many times the following lemma 
\cite[2.3.4 p. 30]{Dixmier}, whose proof is included for the reader's
convenience:

\begin{lemma}[Schur]
\label{Schur}Let $\pi _{1}$ and $\pi _{2}$ be two irreducible
representations of a C*-algebra $A$ acting respectively on Hilbert spaces $%
\mathcal{H}_{1}$ and $\mathcal{H}_{2}$. Then $\pi _{1}$ and $\pi _{2}$ are
unitarily equivalent if and only if there exists a nonzero operator $T:%
\mathcal{H}_{2}\longrightarrow \mathcal{H}_{1}$ such that for all $a\in A$
we have $T\pi _{1}(a)=\pi _{2}(a)T$. Moreover, if there exists such a
nonzero intertwining operator, then it is unique up to a nonzero scalar
multiple.
\end{lemma}

\begin{proof}
If $\pi _{1}$ and $\pi _{2}$ are unitarily equivalent then there exists a
unitary $T$ such that for all $a\in A$ we have $T\pi _{1}(a)=\pi _{2}(a)T$.
In particular, $T\not=0$. Moreover, assume that there exists $T^{\prime }$
such that $T^{\prime }\pi _{1}=\pi _{2}T^{\prime }$. Then $T^{\ast
}T^{\prime }\pi _{1}=T^{\ast }\pi _{2}T^{\prime }=\pi _{1}T^{\ast }T^{\prime
}$. Hence since $\pi _{1}$ is irreducible, there exists $\lambda \in \mathbb{%
C}$ such that $T^{\prime }=\lambda T$.

Conversely, assume that there exists a nonzero operator $T:\mathcal{H}%
_{2}\longrightarrow \mathcal{H}_{1}$ such that for all $a\in A$ we have:%
\begin{equation}
T\pi _{1}(a)=\pi _{2}(a)T\text{.}  \label{UnitaryEq}
\end{equation}%
Then for all $a\in A$:%
\begin{equation*}
T^{\ast }T\pi _{1}(a)=T^{\ast }\pi _{2}(a)T\text{.}
\end{equation*}%
In particular $T^{\ast }T\pi _{1}(a^{\ast })=T^{\ast }\pi _{2}(a^{\ast })T$
for all $a\in A$. Applying the adjoint operation to this equality leads to $%
\pi _{1}(a)T^{\ast }T=T^{\ast }\pi _{2}(a)T$ and thus:%
\begin{equation*}
T^{\ast }T\pi _{1}(a)=\pi _{1}(a)T^{\ast }T\text{.}
\end{equation*}%
Since $\pi _{1}$ is irreducible, there exists $\lambda \in \mathbb{C}$ such
that $T^{\ast }T=\lambda 1_{\mathcal{H}_{2}}$. Since $T\not=0$ we have $%
\lambda \not=0$. Up to replacing $T$ by $\frac{1}{\mu }T$ where $\mu
^{2}=\left\vert \lambda \right\vert $ and $\mu \in \mathbb{R}$ we thus get $%
T^{\ast }T=1_{\mathcal{H}_{2}}$. Thus $T$ is an isometry. In particular, $%
TT^{\ast }$ is a nonzero projection.

Similarly, we get $\pi _{2}(a)TT^{\ast }=TT^{\ast }\pi _{2}(a)$ and thus $%
TT^{\ast }$ is scalar as well. Hence $TT^{\ast }$ is the identity again (As
the only nonzero scalar projection) and thus $T$ is a unitary operator.
Hence by (\ref{UnitaryEq}), $\pi _{1}$ and $\pi _{2}$ are unitarily
equivalent.
\end{proof}

\bigskip Given a Hilbert space $\mathcal{H}$, the C*-algebra of all bounded
linear operators on $\mathcal{H}$ is denoted by $\mathcal{B}\left( \mathcal{H%
}\right) $.

\begin{theorem}
\label{RegularIrred}Let $G$ be a finite group with neutral element $e$ and $%
\alpha $ an action of $G$ on a unital C*-algebra $A$. Let $\pi :A\rightarrow 
\mathcal{B}\left( \mathcal{H}\right) $ be a representation of $A$ and let $%
\lambda $ be the left regular representation of $G$ on $\ell _{2}(G)$. Let $%
\delta _{g}$ be the function in $\ell _{2}(G)$ which is $1$ at $g\in G$ and $%
0$ otherwise. Define $\Pi :A\rtimes _{\alpha }G\rightarrow \mathcal{B}\left(
\ell _{2}\left( G\right) \otimes \mathcal{H}\right) $ by%
\begin{eqnarray*}
\Pi \left( a\right) \left( \delta _{g}\otimes \xi \right) &=&\delta
_{g}\otimes \pi \left( \alpha _{g^{-1}}\left( a\right) \right) \xi ,\text{
and} \\
\Pi \left( g\right) &=&\lambda \left( g\right) \otimes 1_{\mathcal{H}}\text{.%
}
\end{eqnarray*}%
Then $\Pi $ is irreducible if and only if $\pi $ is irreducible and $\pi $
is not unitarily equivalent to $\pi \circ \alpha _{g}$ for any $g\in
G\setminus \left\{ e\right\} $.
\end{theorem}

\begin{proof}
Assume now that $\pi $ is irreducible and not unitarily equivalent to $\pi
\circ \alpha _{g}$ whenever $g\in G\setminus \left\{ e\right\} $. Suppose
that $\Pi $ is reducible. Then there exists a non-scalar operator $\Omega $
in the commutant of $\Pi \left( A\rtimes _{\alpha }G\right) $. Now, we
observe that the commutant of $\left\{ \lambda \left( g\right) \otimes 1_{%
\mathcal{H}}:g\in G\right\} $ is $\rho \left( G\right) \otimes \mathcal{B}%
\left( \mathcal{H}\right) $, where $\rho $ is the right regular
representation of $G$. Hence, there exist an operator $T_{g}$ on $\mathcal{H}
$ for all $g\in G$ such that $\Omega =\sum_{g\in G}$ $\rho \left( g\right)
\otimes T_{g}$. For every $\xi \in \mathcal{H}$ and $a\in A$, we have%
\begin{eqnarray*}
\left( \sum_{g\in G}\rho \left( g\right) \otimes T_{g}\right) \Pi \left(
a\right) \left( \delta _{0}\otimes \xi \right) &=&\left( \sum_{g\in G}\rho
\left( g\right) \otimes T_{g}\right) \left( \delta _{0}\otimes \pi \left(
a\right) \xi \right) \\
&=&\sum_{g\in G}\delta _{g}\otimes T_{g}\pi \left( a\right) \xi
\end{eqnarray*}

and%
\begin{eqnarray*}
\Pi \left( a\right) \left( \sum_{g\in G}\rho \left( g\right) \otimes
T_{g}\right) \left( \delta _{0}\otimes \xi \right) &=&\Pi \left( a\right)
\left( \sum_{g\in G}\delta _{g}\otimes T_{g}\xi \right) \\
&=&\sum_{g\in G}\delta _{g}\otimes \pi \left( \alpha _{g^{-1}}(a)\right)
T_{g}\xi \text{.}
\end{eqnarray*}

Therefore, for every $g\in G$ and for all $a\in A$:%
\begin{equation}
\pi \left( \alpha _{g^{-1}}(a)\right) T_{g}=T_{g}\pi \left( a\right) \text{.}
\label{RegularIrred1}
\end{equation}%
\qquad

Since $\Omega $ is non scalar, there exists $g_{0}\in G\setminus \left\{
e\right\} $ such that $T_{g_{0}}\not=0$. By Lemma (\ref{Schur}), Equality (%
\ref{RegularIrred1}) for $g_{0}$ implies that $\pi $ and $\pi \circ \alpha
_{g_{0}}$, which are irreducible, are also unitarily equivalent since $%
T_{g_{0}}\not=0$. This is a contradiction. So $\Pi $ is irreducible.

We now show the converse. First, note that if $\pi $ is reducible then there
exists a projection $p$ on $\mathcal{H}$ which is neither $0$ or $1$ such
that $p$ commutes with the range of $\pi $. It is then immediate that $%
1\otimes p$ commutes with the range of $\Pi $ and thus $\Pi $ is reducible.

Assume now that there exists $g\in G\setminus \left\{ e\right\} $ such that $%
\pi $ and $\pi \circ \alpha _{g}$ are unitarily equivalent. Then there
exists a unitary $V$ such that for every $a\in A$: 
\begin{equation*}
\pi \left( a\right) =V\pi \left( \alpha _{g}\left( a\right) \right) V^{\ast }%
\text{.}
\end{equation*}%
Let us show that $\rho \left( g\right) \otimes V$ is in the commutant of $%
\Pi \left( A\rtimes _{\alpha }G\right) .$ We only need to check that it
commutes with $\Pi \left( a\right) $ for $a\in A$.%
\begin{eqnarray*}
\left( \rho \left( g\right) \otimes V\right) \Pi \left( a\right) \left(
\delta _{h}\otimes \xi \right) &=&\delta _{hg}\otimes V\pi \left( \alpha
_{h^{-1}}\left( a\right) \right) \xi \text{, and} \\
\Pi \left( a\right) \left( \rho \left( g\right) \otimes V\right) \left(
\delta _{h}\otimes \xi \right) &=&\delta _{hg}\otimes \pi \left( \alpha
_{g^{-1}}\alpha _{h^{-1}}\left( a\right) \right) V\xi \text{.}
\end{eqnarray*}%
Since $V\pi \left( \alpha _{h^{-1}}\left( a\right) \right) =\pi \left(
\alpha _{g^{-1}}\alpha _{h^{-1}}\left( a\right) \right) V,$ we conclude that
the two quantities are equal, and that $\Pi $ is reducible.

Hence, if $\Pi $ is irreducible, then $\pi $ is irreducible and not
equivalent to $\pi \circ \alpha _{g}$ for any $g\in G\setminus \left\{
e\right\} $.
\end{proof}

\bigskip Theorem (\ref{RegularIrred}) provides us with a possible family of
irreducible representations of the crossed-product. The representations
given in Theorem (\ref{RegularIrred}) are called \emph{regular
representations} of $A\rtimes _{\alpha }G$, whether or not they are
irreducible.

\bigskip However, we shall see that there are many irreducible
representations of $A\rtimes _{\alpha }G$ which are not regular. Easy
examples are provided by actions of finite cyclic groups by inner
automorphisms on full matrix algebras, where the identity representation is
in fact the only irreducible representation of the crossed-product. More
generally, the conditions that $\pi $ is irreducible and $\pi \circ \alpha
_{g}$ are not equivalent for $g\in G\backslash \{e\}$ are not necessary.
These observations will be placed into a more general context as we now
address the question raised at the start of this paper in the next section.

\section{Actions of Finite Groups}

This section is concerned with establishing results describing the
irreducible representations of crossed-products by finite groups. The main
tool for our study is to understand such actions from the perspective of the
spectrum of the C*-algebra. In this paper, the spectrum $\widehat{A}$ of a
C*-algebra $A$ is the set of unitary equivalence classes of irreducible
representations of $A$.

We start by two simple observations. Let $\alpha $ be the action of a finite
group $G$ on some unital C*-algebra $A$. Let $\pi _{1}$ and $\pi _{2}$ be
two equivalent irreducible representations of $A$, so that there exists a
unitary $u$ such that $u\pi _{1}u^{\ast }=\pi _{2}$. Then trivially $u\left(
\pi _{1}\circ \alpha _{g^{-1}}\right) u^{\ast }=\pi _{2}\circ \alpha
_{g^{-1}}$ for all $g\in G$. Moreover, $\pi _{1}\circ \alpha _{g^{-1}}$ has
the same range as $\pi _{1}$ and thus is irreducible as well. These two
remarks show that for all $g\in G$ there exists a map $\widehat{\alpha _{g}}$
of $G$ on $\widehat{A}$ defined by mapping the class of an irreducible
representation $\pi $ of $A$ to the class of $\pi \circ \alpha _{g^{-1}}$.
Since $\left( \pi \circ \alpha _{g^{-1}}\right) \circ \alpha _{h^{-1}}=\pi
\circ \alpha _{\left( hg\right) ^{-1}}$, we have $\widehat{\alpha }_{h}\circ 
\widehat{\alpha }_{g}=\widehat{\alpha }_{hg}$, and trivially $\widehat{%
\alpha }_{e}$ is the identity on $\widehat{A}$. Thus $\widehat{\alpha }$ is
an action of $G$ on $\widehat{A}$.

Given a representation $\Pi $ of the crossed-product $A\rtimes _{\alpha }G$,
we define the support of $\Pi $ as the subset $\Sigma $ of $\widehat{A}$ of
all classes of irreducible representations of $A$ weakly contained in $\Pi
_{|A}$. Our main interest are in the support of irreducible representations
of $A\rtimes _{\alpha }G$ which we now prove are always finite.

\subsection{Finiteness of irreducible supports}

Let $G$ be a finite group of neutral element $e$. Let $\widehat{G}$ be the
dual of $G$ i.e. the set of unitary equivalence classes of irreducible
representations of $G$. By \cite[15.4.1, p. 291]{Dixmier}, the cardinal of $%
\widehat{G}$ is given by the number of conjugacy classes of $G$, so $%
\widehat{G}$ is a finite set. Let $\rho \in \widehat{G}$ and $\lambda $ be
any irreducible representation of $G$ of class $\rho $ acting on a Hilbert
space $\mathcal{H}$. Then $\overline{\lambda }$ is the (irreducible)\
representation $g\in G\mapsto \lambda (g)$ acting on the conjugate Hilbert
space $\overline{\mathcal{H}}$ \cite[13.1.5, p. 250]{Dixmier}. We define $%
\overline{\rho }$ as the class of representations unitarily equivalent to $%
\overline{\lambda }$.

\bigskip Let $B$ be a unital C*-algebra and $\alpha $ an action of $G$ on $B$
by *-automorphisms. We now recall from \cite{Hoegh-Krohn81} the definition
and elementary properties of the spectral subspaces of $B$ for the action $%
\alpha $ of $G$. Let $\rho \in \widehat{G}$. The character of $\rho $ is
denoted by $\chi _{\rho }$. All irreducible representations of $G$ whose
class in $\widehat{G}$ is $\rho $ act on vector spaces of the same dimension
which we denote by $\dim \rho $. We recall from \cite[15.3.3, p. 287]%
{Dixmier} that for any $\rho ,\rho ^{\prime }\in \widehat{G}$ we have:%
\begin{equation*}
\chi _{\rho }(e)=\dim \rho
\end{equation*}%
and:%
\begin{equation*}
\chi _{\rho }\ast \chi _{\rho ^{\prime }}(g)=\sum_{h\in G}\chi _{\rho
}(h)\chi _{\rho ^{\prime }}(gh^{-1})=\left\{ 
\begin{array}{ccc}
0 & \text{if} & \rho \not=\rho ^{\prime }\text{,} \\ 
\left( \dim \rho \right) ^{-1}\chi _{\rho }(g) & \text{if} & \rho =\rho
^{\prime }\text{.}%
\end{array}%
\right.
\end{equation*}

The spectral subspace of $B$ for $\alpha $ associated to $\rho \in \widehat{G%
}$ is the space $B_{\rho }$ defined by:%
\begin{equation*}
B_{\rho }=\left\{ \frac{\dim \left( \rho \right) }{\left\vert G\right\vert }%
\sum_{g\in G}\chi _{\overline{\rho }}(g)\alpha _{g}(b):b\in B\right\} \text{,%
}
\end{equation*}%
i.e. the range of the Banach space operator on $B$ defined by: 
\begin{equation}
P_{\rho }:b\in B\mapsto \frac{\dim \left( \rho \right) }{\left\vert
G\right\vert }\sum_{g\in G}\chi _{\overline{\rho }}(g)\alpha _{g}(b)\text{.}
\label{SpectralProjDef}
\end{equation}%
In particular, the spectral subspace associated to the trivial
representation is the fixed point C*-subalgebra $B_{1}$ of $B$ for the
action $\alpha $ of $G$. Now, we have:%
\begin{eqnarray}
P_{\rho }\left( P_{\rho ^{\prime }}(a)\right) &=&\frac{\dim \left( \rho
\right) }{\left\vert G\right\vert }\frac{\dim \left( \rho ^{\prime }\right) 
}{\left\vert G\right\vert }\sum_{g\in G}\sum_{h\in G}\chi _{\overline{\rho }%
}(g)\chi _{\overline{\rho ^{\prime }}}(h)\alpha _{gh}(a)  \notag \\
&=&\frac{\dim \left( \rho \right) }{\left\vert G\right\vert }\frac{\dim
\left( \rho ^{\prime }\right) }{\left\vert G\right\vert }\sum_{g\in G}\left(
\sum_{h\in G}\chi _{\overline{\rho }}(gh^{-1})\chi _{\overline{\rho ^{\prime
}}}(h)\right) \alpha _{g}(a)  \notag \\
&=&\left\{ 
\begin{array}{ccc}
0 & \text{if} & \rho \not=\rho ^{\prime } \\ 
\frac{\dim \left( \rho \right) }{\left\vert G\right\vert }\sum_{g\in G}\chi
_{\overline{\rho }}(g)\alpha _{g}(a) & \text{if} & \rho =\rho ^{\prime }%
\text{.}%
\end{array}%
\right.  \label{SpectralOrthogonal}
\end{eqnarray}%
Hence $P_{\rho }^{2}=P_{\rho }$ so $P_{\rho }$ is a Banach space projection
and $P_{\rho }P_{\rho ^{\prime }}=0$ for all $\rho ^{\prime }\not=\rho $ so
these projections are pairwise orthogonal.

Moreover, for any $g,h\in G$, from \cite[15.4.2 (2) p. 292]{Dixmier}:%
\begin{equation}
\sum_{\rho \in \widehat{G}}\chi _{\rho }(g)\overline{\chi _{\rho }(h)}%
=\left\{ 
\begin{array}{cc}
\frac{\left\vert G\right\vert }{C(g)} & \text{if }g\text{ is conjugated with 
}h\text{,} \\ 
0 & \text{otherwise.}%
\end{array}%
\right.  \label{Sumation2}
\end{equation}%
where for $g\in G$ the quantity $C(g)$ is the number of elements in $G$
conjugated to $g$. In particular, note that since $g\in G\backslash \left\{
e\right\} $ is not conjugated to $e$, we have by Equality (\ref{Sumation2})\
that:%
\begin{equation}
\sum_{\rho \in \widehat{G}}\chi _{\rho }(g)\dim \rho =\sum_{\rho \in 
\widehat{G}}\chi _{\rho }(g)\overline{\chi _{\rho }(e)}=0\text{.}
\label{NullSum}
\end{equation}

Furthermore, because each irreducible representation $\rho $ of $G$ appears
with multiplicity $\dim \rho $ in the left regular representation of $G$ one
can show \cite[15.4.1, p. 291]{Dixmier} that:%
\begin{equation}
\sum_{\rho \in \widehat{G}}\left( \dim \rho \right) ^{2}=\left\vert
G\right\vert \text{.}  \label{RegularMul}
\end{equation}

Hence for all $b\in B$:%
\begin{eqnarray}
\sum_{\rho \in \widehat{G}}P_{\rho }(b) &=&\sum_{\rho \in \widehat{G}}\frac{%
\dim \left( \rho \right) }{\left\vert G\right\vert }\sum_{g\in G}\chi _{\rho
}(g)\alpha _{g}(b)  \notag \\
&=&\frac{1}{\left\vert G\right\vert }\sum_{g\in G}\left( \sum_{\rho \in 
\widehat{G}}\dim \left( \rho \right) \chi _{\rho }(g)\right) \alpha _{g}(b) 
\notag \\
&=&\frac{1}{\left\vert G\right\vert }\left( \sum_{\rho \in \widehat{G}}\dim
\left( \rho \right) \chi (e)\right) \alpha _{e}(b)\text{ by Equality (\ref%
{NullSum})}  \notag \\
&=&\left( \frac{1}{\left\vert G\right\vert }\sum_{\rho \in \widehat{G}}\dim
(\rho )^{2}\right) b=b\text{ by Equality (\ref{RegularMul}).}
\label{ProjectionSummation}
\end{eqnarray}%
Hence $\sum_{\rho \in \widehat{G}}P_{\rho }=\limfunc{Id}_{B}$. Thus by (\ref%
{SpectralOrthogonal})\ and (\ref{ProjectionSummation}) we have:%
\begin{equation}
B=\dbigoplus\limits_{\rho \in \widehat{G}}B_{\rho }\text{.}
\label{SpectralSummation}
\end{equation}

We now establish that the restriction of any irreducible representation of a
crossed-product of some unital C*-algebra $A$ by $G$ is the direct sum of
finitely many irreducible representations of $A$.

\begin{theorem}
\label{FiniteRep}Let $G$ be a finite group and $A$ a unital C*-algebra. Let $%
\alpha $ be an action of $G$ by *-automorphism on $A$. Let $\Pi $ be an
irreducible representation of $A\rtimes _{\alpha }G$ on some Hilbert space $%
\mathcal{H}$. We denote by $U^{g}$ the canonical unitary in $A\rtimes
_{\alpha }G$ corresponding to $g\in G$. Then:

\begin{itemize}
\item The action $g\mapsto \limfunc{Ad}\Pi (U^{g})$ on $\mathcal{B}\left( 
\mathcal{H}\right) $ leaves the commutant $\Pi \left( A\right) ^{\prime }$
of $\Pi (A)$ invariant, and thus defines an action $\beta $ of $G$ on $\Pi
\left( A\right) ^{\prime }$,

\item The action $\beta $ is ergodic on $\Pi \left( A\right) ^{\prime }$,

\item The Von Neumann algebra $\Pi \left( A\right) ^{\prime }$ is finite
dimensional,

\item The representation $\Pi _{|A}$ of $A$ is equivalent to the direct sum
of finitely many irreducible representations of $A$.
\end{itemize}
\end{theorem}

\begin{proof}
Let $\mathfrak{M}=\Pi (A)^{\prime }$. Denote $U_{\Pi }^{g}=\Pi (U^{g})$ for
all $g\in G$. Let $T\in \mathfrak{M}$. Let $a\in A$ and $g\in G$. Then:%
\begin{eqnarray*}
U_{\Pi }^{g}TU_{\Pi }^{g\ast }\Pi (a) &=&U_{\Pi }^{g}TU_{\Pi }^{g\ast }\Pi
(a)U_{\Pi }^{g}U_{\Pi }^{g\ast } \\
&=&U_{\Pi }^{g}T\Pi \left( \alpha _{g^{-1}}(a)\right) U_{\Pi }^{g\ast } \\
&=&U_{\Pi }^{g}\Pi \left( \alpha _{g^{-1}}(a)\right) TU_{\Pi }^{g\ast } \\
&=&U_{\Pi }^{g}U_{\Pi }^{g\ast }\Pi (a)U_{\Pi }^{g}TU_{\Pi }^{g\ast } \\
&=&\Pi (a)U_{\Pi }^{g}TU_{\Pi }^{g\ast }\text{.}
\end{eqnarray*}%
Hence $U_{\Pi }^{g}TU_{\Pi }^{g\ast }\in \mathfrak{M}$ for all $g\in G$ and $%
T\in \mathfrak{M}$. Define $\beta _{g}(T)=U_{\Pi }^{g}TU_{\Pi }^{g\ast }$
for all $g\in G$ and $T\in \mathfrak{M}$. Then $g\in G\mapsto \beta _{g}$ is
an action of $G$ on $\mathfrak{M}$.

Let now $T\in \mathfrak{M}$ such that $\beta _{g}(T)=T$ for all $g\in G$.
Then $T$ commutes with $U_{\Pi }^{g}$ for all $g\in G$. Moreover by
definition of $\mathfrak{M}$, the operator $T$ commutes with $\Pi (A)$.
Hence $T$ commutes with $\Pi $ which is irreducible, so $T$ is scalar. Hence 
$\beta $ is ergodic.

Let $\rho $ be an irreducible representation of $G$ (since $G$ is finite, $%
\rho $ is finite dimensional). By \cite[Proposition 2.1]{Hoegh-Krohn81}, the
spectral subspace $\mathfrak{M}_{\rho }$ of $\mathfrak{M}$ for $\beta $
associated to $\rho $ is finite dimensional. Since $\mathfrak{M}=\oplus
_{\rho \in \widehat{G}}\mathfrak{M}_{\rho }$ by Equality (\ref%
{SpectralSummation}) and since $\widehat{G}$ is finite by \cite[15.4.1, p.
291]{Dixmier} we conclude that $\mathfrak{M}$ is finite dimensional.

Denote $\Pi _{|A}$ by $\pi _{A}$. Let $p_{1},\ldots ,p_{k}$ be projections
in $\mathfrak{M}$, all minimal and such that $\sum_{i=1}^{k}p_{i}=1$. Let $%
i\in \left\{ 1,\ldots ,k\right\} $. Then by definition of $\mathfrak{M}$,
the projection $p_{i}$ commutes with $\pi _{A}$. Hence $p_{i}\pi _{A}p_{i}$
is a representation of $A$. Let $q$ be a projection of $p_{i}\mathcal{H}$
such that $p_{i}$ commutes with $p_{i}\pi _{A}p_{i}$. Then $q\leq p_{i}$ and 
$q\in \mathfrak{M}$, so $q\in \left\{ 0,p_{i}\right\} $ since $p_{i}$ is
minimal. Hence $p_{i}\pi _{A}p_{i}$ is an irreducible representation of $A$.
Therefore:%
\begin{eqnarray*}
\pi _{A} &=&\left( \sum_{i=1}^{k}p_{i}\right) \pi _{A}\text{ since }%
\sum_{i=1}^{k}p_{i}=1\text{,} \\
&=&\sum_{i=1}^{k}p_{i}\pi _{A}p_{i}\text{ since }p_{i}=p_{i}^{2}\in 
\mathfrak{M}\text{.}
\end{eqnarray*}%
Hence $\pi _{A}$ is the direct sum of finitely many irreducible
representations of $A$.
\end{proof}

\subsection{Minimality of the irreducible supports}

\bigskip The following is our key observation which will drive the proofs in
this section:

\begin{Observation}
\label{Observation}Let $\Pi $ be an irreducible representation of $A\rtimes
_{\alpha }G$ and let $\pi _{A}=\Pi _{|A}$. Then for each $g\in G$ the
representations $\pi _{A}$ and $\pi _{A}\circ \alpha _{g}$ are unitarily
equivalent. Hence, the decompositions in direct sums of irreducible
representations of $A$ for $\pi _{A}$ and $\pi _{A}\circ \alpha _{g}$ are
the same.
\end{Observation}

\bigskip This observation is the basis of the next lemma, which is
instrumental in the proof of the theorem to follow.

\begin{lemma}
\label{LemmaCycle}Let $\alpha $ be an action of a finite group $G$ on a
unital C*-algebra $A$. Let $\Pi $ be an irreducible representation of $%
A\rtimes _{\alpha }G$ and let $\pi _{A}$ be the restriction of $\Pi $ to $A$%
. Then there exists a finite subset $\Sigma $ of the spectrum $\widehat{A}$
of $A$ such that all irreducible subrepresentations of $\pi _{A}$ are in $%
\Sigma $. Moreover, all the elements of $\Sigma $ in a given orbit for $%
\widehat{\alpha }$ have the same multiplicity in $\pi _{A}$.
\end{lemma}

\begin{proof}
Let $\Sigma $ be the subset of the spectrum $\widehat{A}$ of $A$ consisting
of all classes of irreducible representations weakly contained in $\pi _{A}$%
. By Theorem (\ref{FiniteRep}), since $\Pi $ is irreducible, $\pi _{A}$ is a
finite direct sum of irreducible representations of $A$ so $\Sigma $ is
nonempty and finite.

Let $g\in G$. Now, by Observation (\ref{Observation}), since $\pi _{A}\circ
\alpha _{g^{-1}}$ is unitarily equivalent to $\pi _{A}$, its decomposition
in irreducible representations is the same as the one for $\pi _{A}$. Thus,
if $\eta \in \Sigma $ then $\widehat{\alpha }_{g}\left( \eta \right) \in
\Sigma $. Since $\widehat{\alpha }_{g}$ is a bijection on $\widehat{A}$ and
thus is injective, and since $\Sigma $ is finite, $\widehat{\alpha }_{g}$ is
a permutation of $\Sigma $.

Let $\Sigma _{\alpha }$ be the orbit of $\varphi \in \Sigma $ under $%
\widehat{\alpha }$ and write $\pi _{A}=\pi _{1}\oplus \ldots \oplus \pi _{k}$
using Theorem (\ref{FiniteRep}), where $\pi _{1},\ldots ,\pi _{k}$ are
irreducible representations of $A$, with the class of $\pi _{1}$ being $%
\varphi $. Now, for $g\in G$, let $n_{1,g},\ldots ,n_{m(g),g}$ be the
integers between $1$ and $k$ such that $\pi _{n_{i,g}}$ is equivalent to $%
\pi _{1}\circ \alpha _{g}$. In particular, $m(g)$ is the multiplicity of $%
\pi _{1}\circ \alpha _{g}$ in $\pi _{A}$. Then $\left( \pi _{n_{1,e}}\oplus
\ldots \oplus \pi _{n_{m(1),e}}\right) \circ \alpha _{g}$ must be the
subrepresentation $\pi _{n_{1,g}}\oplus \ldots \oplus \pi _{n_{m(g),g}}$ of $%
\pi _{A}$. So $m(g)=m(e)$ by uniqueness of the decomposition. Hence for all $%
g$ the multiplicity of $\widehat{\alpha }_{g}(\varphi )$ is the same as the
multiplicity of $\varphi $.
\end{proof}

\bigskip We now establish the main theorem of this paper, describing the
structure of irreducible representations of crossed-products by finite
groups. A \emph{unitary projective representation }of $G$ is a map $\Lambda $
from $G$ into the group of unitaries on some Hilbert space such that there
exists a complex valued 2-cocycle $\sigma $ on $G$ satisfying for all $%
g,h\in G$ the identity $\Lambda _{gh}=\sigma (g,h)\Lambda _{g}\Lambda _{h}$.

\begin{theorem}
\label{FiniteGroupConclusion}Let $G$ be a finite group and $\alpha $ be an
action of $G$ on a unital C*-algebra $A$ by *-automorphisms. Let $\Pi $ be
an irreducible representation of $A\rtimes _{\alpha }G$ on some Hilbert
space $\mathcal{H}$. Then there exists a subgroup $H$ of $G$ and a
representation $\pi $ of $A$ on some Hilbert space $\mathcal{J}$ such that,
up to conjugating $\Pi $ by some fixed unitary, and denoting the index of $H$
in $G$ by $m=G:H$ we have the following:

For any subset $\left\{ g_{1},\ldots ,g_{m}\right\} $ of $G$ such that $%
g_{1} $ is the neutral element of $G$ and $Hg_{j}\cap Hg_{i}=\left\{
g_{1}\right\} $ for $i\not=j$ while $G=\cup _{j=1}^{m}Hg_{j}$, we have:

\begin{enumerate}
\item The representations $\pi \circ \alpha _{g_{i}}$ and $\pi \circ \alpha
_{g_{j}}$ are disjoint for $i,j\in \left\{ 1,\ldots ,m\right\} $ and $%
i\not=j $ (so in particular, they are not unitarily equivalent),

\item There exists an irreducible representation $\pi _{1}$ of $A$ on a
Hilbert subspace $\mathcal{H}_{1}$ of $\mathcal{J}$ and some integer $r$
such that $\mathcal{J}=\mathbb{C}^{r}\otimes \mathcal{H}_{1}$ and $\pi =1_{ 
\mathbb{C}^{r}}\otimes \pi _{1}$,

\item For any $h\in H$ there exists a unitary $V^{h}$ on $\mathcal{H}_{1}$
such that $V^{h}\pi _{1}\left( V^{h}\right) ^{\ast }=\pi _{1}\circ \alpha
_{h}$, and $h\in H\mapsto V^{h}$ is a unitary projective representation of $%
H $ on $\mathcal{H}_{1}$,

\item We have $\mathcal{H}=\mathcal{J}_{g_{1}}\oplus \ldots \oplus \mathcal{J%
}_{g_{m}}$ where for all $i=1,\ldots ,m$ the space $\mathcal{J}_{g_{i}}$ is
an isometric copy of $\mathcal{J}$,

\item In this decomposition of $\mathcal{H}$ we have for all $a\in A$ that:%
\begin{equation}
\Pi (a)=\left[ 
\begin{array}{cccc}
\pi (a) &  &  &  \\ 
& \pi \circ \alpha _{g_{2}}(a) &  &  \\ 
&  & \ddots &  \\ 
&  &  & \pi \circ \alpha _{g_{m}}(a)%
\end{array}%
\right]  \label{PIA}
\end{equation}

\item In this same decomposition, for every $g$ there exists a permutation $%
\sigma ^{g}$ of $\left\{ 1,\ldots ,m\right\} $ and unitaries $U_{i}^{g}:%
\mathcal{J}_{g_{i}}\longrightarrow \mathcal{J}_{\sigma ^{g}(g_{i})}$ such
that:%
\begin{equation*}
\Pi (U^{g})=\left[ U_{j}^{g}\delta _{i}^{\sigma ^{g}(j)}\right]
_{i,j=1,\ldots ,m}
\end{equation*}%
where $\delta $ is the Kronecker symbol:%
\begin{equation}
\delta _{a}^{b}=\left\{ 
\begin{array}{cc}
1 & \text{if }a=b\text{,} \\ 
0 & \text{otherwise.}%
\end{array}%
\right.  \label{Kronecker}
\end{equation}%
Moreover:%
\begin{equation*}
H=\left\{ g\in G:\sigma ^{g}(1)=1\right\} \text{.}
\end{equation*}

\item The representation $\Psi $ of $A\rtimes _{\alpha }H$ on $\mathcal{J}$
defined by $\Psi (a)=\pi (a)$ for all $a\in A$ and $\Psi (U^{h})=U_{1}^{h}$
for $h\in H$ is irreducible. Moreover, there exists an irreducible unitary
projective representation $\Lambda $ of $G$ on $\mathbb{C}^{r}$ such that on 
$\mathcal{J}=\mathbb{C}^{r}\otimes \mathcal{H}_{1}$, while $\Psi (a)=1_{%
\mathbb{C}^{r}}\otimes \pi _{1}(a)$, we also have $\Psi
(U^{h})=U_{1}^{h}=\Lambda _{h}\otimes V^{h}$.
\end{enumerate}
\end{theorem}

\begin{proof}
Let $\Pi $ be an irreducible representation of $A\rtimes _{\alpha }G$.
Denote $\Pi _{|A}$ by $\pi _{A}$. By Theorem (\ref{FiniteRep}), there exists
a nonzero natural integer $k$ and irreducible representations $\pi
_{1},\ldots ,\pi _{k}$ of $A$, acting respectively on Hilbert spaces $%
\mathcal{H}_{1},\ldots ,\mathcal{H}_{k}$ such that up to a unitary
conjugation of $\Pi $, we have $\mathcal{H}=\mathcal{H}_{1}\oplus \ldots
\oplus \mathcal{H}_{k}$ and in this decomposition, for all $a\in A$:%
\begin{equation*}
\pi _{A}(a)=\left[ 
\begin{array}{cccc}
\pi _{1}(a) &  &  &  \\ 
& \pi _{2}(a) &  &  \\ 
&  & \ddots &  \\ 
&  &  & \pi _{k}(a)%
\end{array}%
\right] \text{.}
\end{equation*}

At this stage, the indexing of the irreducible subrepresentations of $\pi
_{A}$ is only defined up to a permutation of $\left\{ 1,\ldots ,k\right\} $.
We start our proof by making a careful choice of such an indexing. To do so,
first choose $\pi _{1}$ arbitrarily among all irreducible subrepresentations
of $\pi _{A}$. Our next step is to set:%
\begin{equation*}
H=\left\{ g\in G:\pi _{1}\circ \alpha _{g}\text{ is equivalent to }\pi
_{1}\right\} \text{.}
\end{equation*}%
We now show that $H$ is a subgroup of $G$. For all $h\in H$ we denote by $%
V^{h}$ the (unique, up to a scalar multiple) unitary such that $V^{h}\pi
_{1}\left( V^{h}\right) ^{\ast }=\pi _{1}\circ \alpha _{h}$. Then if $g,h\in
H$ we have:%
\begin{eqnarray*}
\pi _{1}\circ \alpha _{gh^{-1}} &=&\left( \pi _{1}\circ \alpha _{g}\right)
\circ \alpha _{h^{-1}}=V^{g}\left( \pi _{1}\circ \alpha _{h^{-1}}\right)
V^{g^{-1}} \\
&=&V^{g}V^{h^{-1}}\pi _{1}V^{h}V^{g^{-1}}
\end{eqnarray*}%
so $\pi _{1}\circ \alpha _{gh^{-1}}$ is unitarily equivalent to $\pi _{1}$
and thus $gh^{-1}\in H$ by definition. Since $H$ trivially contains the
neutral element of $G$, we conclude that $H$ is a subgroup of $G$.

Let $\left\{ g_{1},\ldots ,g_{m}\right\} $ a family of right coset
representatives such that $g_{1}$ is the neutral element of $G$ \cite[p. 10]%
{Robinson82}, i.e. such that for $i\not=j$ we have $Hg_{j}\cap
Hg_{i}=\left\{ g_{1}\right\} $ while $G=\cup _{j=1}^{m}Hg_{j}$. In
particular, for $i\in \left\{ 2,\ldots ,m\right\} $ we have $g_{1}\not=g_{i}$
and by definition of $H$ this implies that $\pi _{1}\circ \alpha _{g_{i}}$
is not equivalent to $\pi _{1}$.

Then let $\pi _{2},\ldots ,\pi _{n_{1}}$ be all the representations
equivalent to $\pi _{1}$. We then choose $\pi _{n_{1}+1}$ to be a
subrepresentation of $\pi _{A}$ equivalent to $\pi _{1}\circ \alpha _{g_{1}}$%
. Again, we let $\pi _{n_{1}+1},\ldots ,\pi _{n_{2}}$ be all the
representations which are equivalent to $\pi _{n_{1}+1}$. More generally, we
let $\pi _{n_{j}+1},\ldots ,\pi _{n_{j+1}}$ be all the subrepresentations of 
$\pi _{A}$ equivalent to $\pi _{1}\circ \alpha _{g_{j}}$ for all $j\in
\{1,\ldots ,m\}$. All other irreducible subrepresentations of $\pi _{A}$
left, if any, are indexed from $n_{m}+1$ to $k$ and we denote their direct
sum by $\Lambda $.

Note that $\Lambda $ contains no subrepresentation equivalent to any
representation $\pi _{1}\circ \alpha _{g}$ for any $g\in G$. Indeed, if $%
g\in G$ then there exists $h\in H$ and a unique $j\in \left\{ 1,\ldots
,m\right\} $ such that $g=hg_{j}$. Thus:%
\begin{equation*}
\pi _{1}\circ \alpha _{g}=\pi _{1}\circ \alpha _{h}\circ \alpha
_{g_{j}}=V^{h}\left( \pi _{1}\circ \alpha _{g_{j}}\right) V^{-h}
\end{equation*}%
and thus $\pi _{1}\circ \alpha _{g}$ is equivalent to one of the
representations $\pi _{1},\ldots ,\pi _{n_{m}}$ by construction. Also note
that if $\pi _{1}\circ \alpha _{g_{i}}$ is equivalent to $\pi _{1}\circ
\alpha _{g_{j}}$ then $g_{i}g_{j}^{-1}\in H$ which contradicts our choice of 
$\left\{ g_{1},\ldots ,g_{m}\right\} $ unless $i=j$. Hence, for $i\not=j$
the representations $\pi _{1}\circ \alpha _{g_{j}}$ and $\pi _{1}\circ
\alpha _{g_{i}}$ are not equivalent.

Now, if $\varphi _{1},\ldots ,\varphi _{m}$ represent the
unitary-equivalence classes of the representations $\pi _{1},\pi _{1}\circ
\alpha _{g_{1}},\ldots ,\pi _{1}\circ \alpha _{g_{m}}$ then $\Sigma
_{1}=\left\{ \varphi _{1},\ldots ,\varphi _{m}\right\} $ is the orbit of $%
\varphi _{1}$ for the action $\widehat{\alpha }$ of $G$ on $\widehat{A}$.
Therefore, there exists $r\geq 1$ such that $n_{j}=jr+1$ for all $j=1,\ldots
,m$ by Lemma (\ref{LemmaCycle}), i.e. all the representations $\pi _{1}\circ
\alpha _{g_{i}}$ ($i=1,\ldots ,m$) have multiplicity $r$ in $\pi _{A}$.

Thus, (up to equivalence on $\Pi $) and writing $\mathcal{H}=\oplus
_{i=1}^{k}\mathcal{H}_{i}$ and in this decomposition:%
\begin{eqnarray}
\pi _{A} &=&\underset{\text{each equivalent to }\pi _{1}}{\underbrace{\pi
_{1}\oplus \ldots \oplus \pi _{r}}}\oplus \underset{\text{each equivalent to 
}\pi _{1}\circ \alpha _{g_{1}}}{\underbrace{\pi _{r+1}\oplus \ldots \oplus
\pi _{2r}}}\oplus \cdots  \label{DecompositionLambda} \\
&&\ldots \oplus \pi _{mr}\oplus \underset{\Lambda }{\underbrace{\pi
_{mr+1}\oplus \ldots \oplus \pi _{k}}}  \notag \\
&=&\underset{\text{disjoint from }\Lambda \text{.}}{\underbrace{\pi
_{1}\oplus \ldots \oplus \pi _{n_{m}}}}\oplus \Lambda \text{.}  \notag
\end{eqnarray}

Let $g\in G$. Still in the decomposition $\mathcal{H}=\mathcal{H}_{1}\oplus
\ldots \oplus \mathcal{H}_{k}$ with our choice of indexing, let us write:%
\begin{equation*}
\Pi \left( U^{g}\right) =U_{\Pi }^{g}=\left[ 
\begin{array}{cccc}
a_{11}^{g} & a_{12}^{g} & \cdots & a_{1k}^{g} \\ 
a_{21}^{g} & a_{22}^{g} & \cdots & a_{2k}^{g} \\ 
\vdots &  & \ddots & \vdots \\ 
a_{k1}^{g} & a_{k2}^{g} & \cdots & a_{kk}^{g}%
\end{array}%
\right]
\end{equation*}%
for some operators $a_{ij}^{g}$ from $\mathcal{H}_{j}$ to $\mathcal{H}_{i}$
with $i,j=1,\ldots ,k$.

Since $U_{\Pi }^{g}\pi _{A}(a)=\pi _{A}(\alpha _{g}(a))U_{\Pi }^{g}$, we can
write:%
\begin{eqnarray}
&&\left[ 
\begin{array}{cccc}
a_{11}^{g}\pi _{1} & a_{12}^{g}\pi _{2} & \cdots & a_{1k}^{g}\pi _{k} \\ 
a_{21}^{g}\pi _{1} & a_{22}^{g}\pi _{2} & \cdots & a_{2k}^{g}\pi _{k} \\ 
\vdots &  & \ddots & \vdots \\ 
a_{k1}^{g}\pi _{1} & a_{k2}^{g}\pi _{2} & \cdots & a_{kk}^{g}\pi _{k}%
\end{array}%
\right]  \label{MainEquality} \\
&=&\left[ 
\begin{array}{cccc}
\left( \pi _{1}\circ \alpha _{g}\right) a_{11}^{g} & \left( \pi _{1}\circ
\alpha _{g}\right) a_{12}^{g} & \cdots & \left( \pi _{1}\circ \alpha
_{g}\right) a_{1k}^{g} \\ 
\left( \pi _{2}\circ \alpha _{g}\right) a_{21}^{g} & \left( \pi _{2}\circ
\alpha _{g}\right) a_{22}^{g} & \cdots & \left( \pi _{2}\circ \alpha
_{g}\right) a_{2k}^{g} \\ 
\vdots &  & \ddots & \vdots \\ 
\left( \pi _{k}\circ \alpha _{g}\right) a_{k1}^{g} & \left( \pi _{k}\circ
\alpha _{g}\right) a_{k2}^{g} & \cdots & \left( \pi _{k}\circ \alpha
_{g}\right) a_{kk}^{g}%
\end{array}%
\right] \text{.}  \notag
\end{eqnarray}%
As a consequence of Equality (\ref{MainEquality}), we observe that for all $%
i,j\in \left\{ 1,\ldots ,k\right\} $ we have:%
\begin{equation}
a_{ij}^{g}\pi _{j}=\left( \pi _{i}\circ \alpha _{g}\right) a_{ij}^{g}\text{.}
\label{SchurEquality}
\end{equation}

First, let $i>mr$. Then the equivalence class of $\pi _{i}$ is not in the
orbit $\Sigma _{1}$ of $\varphi _{1}$ for $\widehat{\alpha }$ by
construction. Hence $\pi _{i}\circ \alpha _{g}$ is not unitarily equivalent
to $\pi _{1}\circ \alpha _{\gamma }$ for any $\gamma \in G$. On the other
hand, let $j\leq mr$. The representation $\pi _{j}$ is equivalent to $\pi
_{1}\circ \alpha _{g_{l}}$ for some $l\in \left\{ 1,\ldots ,m\right\} $ by
our choice of indexing. Therefore, $\pi _{i}\circ \alpha _{g}$ and $\pi _{j}$
are not unitarily equivalent, yet they both are irreducible representations
of $A$. Hence by Lemma (\ref{Schur}) applied to Equality (\ref{SchurEquality}%
) we conclude that $a_{ij}^{g}=0$. Similarly, $\pi _{i}$ and $\pi _{j}\circ
\alpha _{g}$ are not equivalent so $a_{ji}^{g}=0$ as well.

Hence:%
\begin{equation*}
U_{\Pi }^{g}=\left[ 
\begin{array}{cccccc}
a_{11}^{g} & \cdots & a_{1mr}^{g} & 0 & \cdots & 0 \\ 
\vdots &  & \vdots & \vdots &  & \vdots \\ 
a_{mr1}^{g} & \cdots & a_{mr,mr}^{g} & 0 & \cdots & 0 \\ 
0 & \cdots & 0 & a_{mr+1,mr+1}^{g} & \cdots & a_{mr+1,k}^{g} \\ 
\vdots &  & \vdots & \vdots & \ddots & \vdots \\ 
0 & \cdots & 0 & a_{k,mr+1}^{g} & \cdots & a_{kk}^{g}%
\end{array}%
\right] \text{.}
\end{equation*}%
If we assume that $n_{m}=mr<k$ then for all $g\in G$ the unitary $U_{\Pi
}^{g}$ commutes with the nontrivial projection $\underset{mr\text{ times}}{%
\underbrace{0\oplus \ldots \oplus 0}}\oplus \underset{k-mr\text{ times}}{%
\underbrace{1\oplus \ldots \oplus 1}}$ of $\mathcal{H}$, and so does $\pi
_{A}$. Yet $\Pi $ is irreducible, so this is not possible and thus $n_{m}=k$%
. Thus $\Sigma =\Sigma _{1}$ is an orbit of a single $\varphi \in \widehat{A}
$ for $\widehat{\alpha }$ and there is no $\Lambda $ left in Equality (\ref%
{DecompositionLambda}). In particular, the cardinal of $\Sigma _{1}$ is $m$.

Since by construction $\pi _{jr+z}$ is unitarily equivalent to $\pi
_{1}\circ \alpha _{g_{z}}$ for all $j=0,\ldots ,m-1$ and $z=1,\ldots ,r$,
there exists a unitary $\omega _{jr+z}$ from $\mathcal{H}_{1}$ onto $%
\mathcal{H}_{jr+z}$ such that $\omega _{jr+z}\left( \pi _{1}\circ \alpha
_{g_{z}}\right) \omega _{jr+z}^{\ast }=\pi _{jr+z}$ (note that we can choose 
$\omega _{1}=1$). We define on $\mathcal{H}=\mathcal{H}_{1}\oplus \ldots
\oplus \mathcal{H}_{k}$ the diagonal unitary:%
\begin{equation*}
\Omega =\left[ 
\begin{array}{ccc}
\omega _{1}^{\ast } &  &  \\ 
& \ddots &  \\ 
&  & \omega _{k}^{\ast }%
\end{array}%
\right] \text{.}
\end{equation*}%
Denote by $\limfunc{Ad}\Omega $ is the *-automorphism on the C*-algebra of
bounded operators on $\mathcal{H}$ defined by $T\mapsto \Omega T\Omega
^{\ast }$. Then up to replacing $\Pi $ by $\limfunc{Ad}\Omega \circ \Pi $,
we can assume that $\pi _{jr+z}=\pi _{1}\circ \alpha _{g_{z}}$ for all $j\in
\left\{ 1,\ldots ,m\right\} $ and $z\in \left\{ 1,\ldots ,r\right\} $. Given
an irreducible representation $\eta $ of $A$ and any nonzero natural integer 
$z$ we shall denote by $z\cdot \eta $ the representation $\underset{z\text{
times}}{\underbrace{\eta \oplus \ldots \oplus \eta }}$. Thus, if we set $\pi
=r\cdot \pi _{1}$ we see that $\pi _{A}$ can be written as in Equality (\ref%
{PIA}) with $\pi \circ \alpha _{g_{i}}$ disjoint from $\pi \circ \alpha
_{g_{j}}$ for $i,j\in \left\{ 1,\ldots ,m\right\} $ and $i\not=j$.

Let again $g\in G$. We now use the same type of argument to show that $%
U_{\Pi }^{g}$ is a \textquotedblleft unitary-permutation
shift\textquotedblright . Let $j\in \left\{ 0,\ldots ,m-1\right\} $. Let $%
q\in \left\{ 1,\ldots ,m\right\} $ such that $g_{j}g\in Hg_{q}$ --- by our
choice of $g_{1},\ldots ,g_{m}$ there is a unique such $q$. Let $i\in
\left\{ 0,\ldots ,m-1\right\} \backslash \left\{ q\right\} $ and $z,h\in
\left\{ 1,\ldots ,r\right\} $. By construction, the representation $\left(
r\cdot \pi _{rj+h}\right) \circ \alpha _{g}$ is unitarily equivalent to $%
r\cdot \pi _{rq+h}$ and disjoint from $r\cdot \pi _{ri+z}$. Yet by Equality (%
\ref{SchurEquality}) we have again that:%
\begin{equation*}
a_{ri+z,rj+h}^{g}\pi _{ri+z}=\left( \pi _{rj+h}\circ \alpha _{g}\right)
a_{ri+z,rj+h}^{g}\text{.}
\end{equation*}%
Thus $a_{ri+z,rj+h}^{g}=0$ by Lemma (\ref{Schur}) since $\pi _{ri+z}$ and $%
\pi _{rj+h}\circ \alpha _{g}$ are not equivalent yet irreducible. Thus, if
for all $z\in \left\{ 0,\ldots ,m-1\right\} $ we define the Hilbert subspace 
$\mathcal{J}_{z}=\mathcal{H}_{zr+1}\oplus \ldots \oplus \mathcal{H}_{(z+1)r}$
of $\mathcal{H}$ then we conclude that $U_{\Pi }^{g}\left( \mathcal{J}%
_{j}\right) \subseteq \mathcal{J}_{q}$ and $\mathcal{H}=\mathcal{J}%
_{0}\oplus \ldots \oplus \mathcal{J}_{m-1}$. Moreover, by uniqueness of $q$
we also obtain that:%
\begin{equation}
U_{\Pi }^{g^{-1}}\left( \mathcal{J}_{q}\right) \subseteq \mathcal{J}_{j}
\label{InverseSubset}
\end{equation}%
and thus $U_{\Pi }^{g}\left( \mathcal{J}_{j}\right) =\mathcal{J}_{q}$.
Define $\sigma ^{g}(j)=q$. Then $\sigma ^{g}$ is a surjection of the finite
set $\left\{ 1,\ldots ,m\right\} $ by (\ref{InverseSubset}), so $\sigma ^{g}$
is a permutation. If $\delta $ is defined as in Equality (\ref{Kronecker})
then, if we set $U_{j}^{g}=U_{\Pi |\mathcal{J}_{j}}^{g}$ then:%
\begin{equation}
\left( r\cdot \left( \pi _{j}\circ \alpha _{g}\right) \right)
U_{i}^{g}=U_{i}^{g}\left( r\cdot \pi _{q}\right)
\end{equation}%
and 
\begin{equation*}
U_{\Pi }^{g}=\left[ U_{j}^{g}\delta _{i}^{\sigma ^{g}(j)}\right] _{i,j}
\end{equation*}
for all $i=1,\ldots ,m$.

Since $U_{\Pi }$ is unitary, so are the operators $U_{1},\ldots ,U_{m}$. In
particular, $\mathcal{J}_{j}$ and $\mathcal{J}_{0}$ are isometric Hilbert
spaces for all $j=0,\ldots ,m-1$. Note that $\left( r\cdot \pi _{1}\right)
\circ \alpha _{g_{i}}$ acts on $\mathcal{J}_{i-1}$ for $i=1,\ldots ,m$ by
construction. We now denote $r\cdot \pi _{1}$ by $\pi $ and $\mathcal{J}=%
\mathcal{J}_{0}$.

Now, by construction $\sigma ^{g}(1)=1$ if and only if there exists an
operator $W$ on $\mathcal{J}_{1}\oplus \ldots \oplus \mathcal{J}_{m-1}$ such
that $U_{\Pi }^{g}=U_{1}^{g}\oplus W$, which is equivalent to $U_{1}^{g}\pi
_{1}=\left( \pi _{1}\circ \alpha _{g}\right) U_{1}^{g}$. By construction,
this is possible if and only if $g\in H$.

Let now $h\in H$. Hence $U_{1}^{h}\pi U_{1}^{h^{-1}}=\pi \circ \alpha _{h}$.
If we set $\Psi (a)=\pi (a)$ and $\Psi (U^{h})=U_{1}^{h}$, we thus define a
representation of $A\rtimes _{\alpha }H$ on $\mathcal{J}_{0}$. Let $b\in
A\rtimes _{\alpha }H$. Then there exists $g\in G\mapsto a_{g}$ such that $%
b=\sum_{g\in G}a_{g}U^{g}$. Hence $\Pi (b)=\sum_{g\in G}\pi
_{A}(a_{g})U_{\Pi }^{g}$. Let $Q$ be the projection of $\mathcal{H}$ on $%
\mathcal{J}_{0}$. Then:%
\begin{eqnarray}
Q\Pi (b)Q &=&\sum_{g\in G}\pi (a_{g})QU_{\Pi }^{g}Q  \notag \\
&=&\sum_{h\in H}\pi (a_{h})U_{1}^{h}=\Psi \left( \sum_{h\in
H}a_{h}U^{h}\right) \text{.}  \label{PsiRange}
\end{eqnarray}%
Since $\Pi $ is irreducible, the range of $\Pi $ is $\limfunc{WOT}$ dense by
the double commutant theorem. Hence, since the multiplication on the left
and right by a fixed operator is $\limfunc{WOT}$ continuous, we conclude
that $Q\Pi Q$ is $\limfunc{WOT}$ dense in $\mathcal{B}\left( Q\mathcal{H}%
\right) $. Therefore, by Equality (\ref{PsiRange}), we conclude by the
double commutant Theorem again that $\Psi $ is an irreducible representation
of $A\rtimes _{\alpha }H$.

Last, note that since $\pi _{1}$ is irreducible, if $h,g\in H$ then since:%
\begin{equation*}
V^{h^{-1}}V^{g^{-1}}V^{gh}\pi _{1}=\pi _{1}V^{h^{-1}}V^{g^{-1}}V^{gh}
\end{equation*}%
there exists $\lambda _{g,h}\in \mathbb{T}$ such that $V^{gh}=\lambda
_{g,h}V^{g}V^{h}$. Hence $g\in H\mapsto V^{g}$ is a projective
representation of $H$ on $\mathcal{H}_{1}$. Note that although the unitaries 
$V^{h}$ are only defined up to a scalar, there is no apparent reason why one
could choose $\lambda $ to be the trivial cocycle unless the second
cohomology group of $H$ is trivial. We now note that $\mathcal{J}_{0}=%
\mathcal{J}=\mathbb{C}^{r}\otimes \mathcal{H}_{1}$ by construction. Now, for
all $h\in H$ we set $\upsilon _{h}=1_{\mathbb{C}^{r}}\otimes V^{h}$. Again, $%
\upsilon _{h}$ is a projective representation of $H$. Moreover, for $h\in H$:%
\begin{equation*}
U_{1}^{h}\upsilon _{h}^{\ast }\pi =\pi U_{1}^{h}\upsilon _{h}^{\ast }\text{.}
\end{equation*}%
Since $\pi =r\cdot \pi _{1}$, Lemma (\ref{Schur}) implies that there exist a
unitary $\Lambda _{h}\in M_{r}\left( \mathbb{C}\right) $ such that $%
U_{1}^{h}\upsilon _{h}^{\ast }=\Lambda _{h}\otimes 1_{\mathcal{H}_{1}}$.
Hence $U_{1}^{h}=\Lambda _{h}\otimes V^{h}$. Now, for $h,g\in H$ we have $%
U_{1}^{h}U_{1}^{g}=U_{1}^{hg}$ which implies that:%
\begin{equation*}
\left( \Lambda _{h}\otimes V^{h}\right) \left( \Lambda _{g}\otimes
V^{g}\right) =\Lambda _{h}\Lambda _{g}\otimes V^{h}V^{g}=\Lambda
_{hg}\otimes V^{hg}\text{.}
\end{equation*}%
Hence $h\mapsto \Lambda _{h}$ is a unitary projective representation of $H$
on $\mathbb{C}^{r}$ with cocycle $\overline{\lambda }$. Moreover, if $T$
commutes with the range of $\Lambda $ then $T\otimes 1$ commutes with the
range of $\Psi $, which contradicts the irreducibility of $\Psi $. Hence $%
\Lambda $ is irreducible. This completes the description of the
representation $\Psi $.
\end{proof}

\bigskip For generic groups, the representation $\Psi $ of Theorem (\ref%
{FiniteGroupConclusion}) may not be minimal, i.e. its restriction to $A$ may
be reducible. The simplest way to see this is by consider a finite group $G$
admitting a representation $\Lambda $ on $\mathbb{C}^{n}$ for some $n\in 
\mathbb{N}$. Then $\Lambda $ extends to an irreducible representation $\Pi $
of the crossed-product $\mathbb{C}\rtimes _{\alpha }G$ where $\alpha $ is
the trivial action. Thus, $\Pi _{|\mathbb{C}}$, which decomposes into a
direct sum of irreducible representations of $\mathbb{C}$, must in fact be
the direct sum of $n$ copies of the (unique) identity representation of $%
\mathbb{C}$. Note that in this case $\Pi =\Psi $ using the notations of
Theorem (\ref{FiniteGroupConclusion}). Thus, for any $n\in \mathbb{N}$ one
can find an example where $\Psi $ is irreducible yet not minimal. This
situation will be illustrated with a much less trivial example in Example (%
\ref{ExPermutation1}) where $G$ will be permutation group on three elements.
However, the representation $\Psi $ must be minimal when the group $G$ is
chosen to be a finite cyclic group. We develop the theory for these groups
in the next section.

Because the representation $\Psi $ of Theorem (\ref{FiniteGroupConclusion})
is of central interest in the decomposition of $\Pi $, we establish the
following criterion for irreducibility for such representations. Note that
the next theorem also describes the situation where the commutant of $\Pi $
is a factor.

\begin{theorem}
\label{Homogeneous}Let $H$ be a discrete group. Let $\Psi $ be a
representation of $A\rtimes _{\alpha }H$ on a Hilbert space $\mathcal{H}$
and assume there exists an irreducible representation $\pi _{1}$ of $A$ on a
Hilbert space $\mathcal{H}_{1}$ such that $\mathcal{H}=\mathbb{C}^{r}\otimes 
\mathcal{H}_{1}$, $\pi _{1}\circ \alpha _{h}$ is equivalent to $\pi _{1}$
for all $h\in H$ and $\Psi (a)=1_{\mathbb{C}^{r}}\otimes \pi (a)$ for all $%
a\in A$. Then there exist two unitary projective representations $\Lambda $
and $V$ of $H$ on $\mathbb{C}^{r}$ and $\mathcal{H}_{1}$ respectively such
that $\Psi (U^{h})=\Lambda _{h}\otimes V^{h}$. Moreover, the following are
equivalent:

\begin{enumerate}
\item $\Psi $ is irreducible,

\item The representation $\Lambda $ is irreducible.
\end{enumerate}
\end{theorem}

\begin{proof}
By assumption, for $h\in H$ there exists a unitary $V^{h}$ such that $%
V^{h}\pi _{1}\left( V^{h}\right) ^{\ast }=\pi _{1}\circ \alpha _{h}$ and
this unitary is unique up to a constant by Lemma (\ref{Schur}). From the
last section of the proof of Theorem (\ref{FiniteGroupConclusion}), we get
that $h\in H\mapsto V^{h}$ is a projective representation of $H$ for some
2-cocycle $\lambda $ and, since $\pi _{1}$ is irreducible, there exists a
projective representation $\Lambda $ of $H$ on $\mathbb{C}^{r}$ such that $%
\Psi (U^{h})=\Lambda _{h}\otimes V^{h}$, and moreover if $\Psi $ is
irreducible then so is $\Lambda $.

Suppose now $\Lambda $ is irreducible. Let $T\in \left[ \Psi \left( A\rtimes
_{\alpha }H\right) \right] ^{\prime }.$ Since $T$ commutes with $\Psi \left(
A\right) =1_{\mathbb{C}^{r}}\otimes \pi _{1}\left( A\right) $, it follows
that $T=D\otimes 1_{\mathcal{H}_{1}}$ for some $D\in M_{r}\left( \mathbb{C}%
\right) $. Now $T$ commutes with $\Psi (U^{h})$ for all $h\in H$, so $D$
commutes with $\Lambda _{g}$ for all $g\in H$. Hence $D$ is scalar and $\Psi 
$ is irreducible.
\end{proof}

We also note that the group $H$ is not a priori a normal subgroup of $G$. It
is easy to check that the following two assertions are equivalent:

\begin{enumerate}
\item $H$ is a normal subgroup of $G$,

\item For all $g\in G$, the unitary $U_{\Pi }^{g}$ is block-diagonal in the
decomposition $\mathcal{H}=\mathcal{J}_{0}\oplus \ldots \oplus \mathcal{J}%
_{m-1}$ if and only if $g\in H$.
\end{enumerate}

In particular, when $G$ is Abelian then for $g\in G$ we have $\sigma
^{g}(1)=1$ if and only if $\sigma ^{g}=\limfunc{Id}$.

\bigskip We conclude by observing that the representation $\Psi $ involves
projective representations of $H$. We now offer an example to illustrate
this situation and shows that this phenomenon occurs even when $G$ is
Abelian. We shall see in the next section that finite cyclic groups have the
remarkable property that such unitary projective representations do not
occur.

\begin{example}
Let $p,q$ be two relatively prime integers. Let $\lambda =\exp \left( 2i\pi 
\frac{p}{q}\right) $. Denote by $\mathbb{U}_{q}$ the group of $q^{\text{th}}$
roots of unity in $\mathbb{C}$. Let $\alpha $ be the action of $\mathbb{Z}%
_{q}$ on $C\left( \mathbb{U}_{q}\right) $ defined by $\alpha
_{1}(f)(z)=f\left( \lambda z\right) $. Then the crossed-product $A=C(\mathbb{%
U}_{q})\rtimes _{\alpha }\mathbb{Z}_{q}$ is isomorphic to $M_{q}(\mathbb{C})$%
. The canonical unitary is identified under this isomorphism with:%
\begin{equation*}
U=\left[ 
\begin{array}{cccc}
0 & 1 & 0 & 0 \\ 
\vdots  &  & \ddots  & \vdots  \\ 
0 &  &  & 1 \\ 
1 & 0 & \cdots  & 0%
\end{array}%
\right] 
\end{equation*}%
while the generator $z\in \mathbb{U}_{q}\mapsto z$ of $C\left( \mathbb{U}%
_{q}\right) $ is mapped to:%
\begin{equation*}
V=\left[ 
\begin{array}{cccc}
1 &  &  &  \\ 
& \lambda  &  &  \\ 
&  & \ddots  &  \\ 
&  &  & \lambda ^{q-1}%
\end{array}%
\right] \text{.}
\end{equation*}%
The dual action $\gamma $ of the Abelian group $G=\mathbb{Z}_{q}\times 
\mathbb{Z}_{q}$ on $C(\mathbb{U}_{q})\rtimes _{\alpha }\mathbb{Z}_{q}$ can
thus be described by:%
\begin{equation*}
\gamma _{z,z^{\prime }}\left( U\right) =\exp \left( 2i\pi \frac{pz}{q}%
\right) U\text{ and }\gamma _{z,z^{\prime }}(V)=\exp \left( 2i\pi \frac{%
pz^{\prime }}{q}\right) V
\end{equation*}%
for all $\left( z,z^{\prime }\right) \in G$. Now, for $(z,z^{\prime })\in G$
we set $\Lambda (z,z^{\prime })=\lambda ^{zz^{\prime }}U^{z}V^{z^{\prime }}$%
. Note that $G$ is generated by $\zeta =\left( 1,0\right) $ and $\xi =\left(
0,1\right) $ and $\Lambda (\zeta )=U$ while $\Lambda (\xi )=V$. Since $%
VU=\lambda UV$, the map $\Lambda $ is a unitary projective representation of 
$G$ on $\mathbb{C}^{2}$ associated to the group cohomology class of $\exp
\left( i\pi \sigma \right) $ where $\sigma $ is defined by:%
\begin{equation*}
\sigma (\left( z,z^{\prime }\right) ,\left( y,y^{\prime }\right) )=\frac{p}{q%
}\left( zy^{\prime }-z^{\prime }y\right) \text{.}
\end{equation*}%
Moreover, the dual action is of course an inner action, and more precisely:%
\begin{eqnarray*}
\gamma _{z,z^{\prime }}\left( a\right)  &=&U^{z}V^{z^{\prime
}}aV^{-z^{\prime }}U^{-z} \\
&=&\Lambda (z,z^{\prime })a\Lambda \left( z,z^{\prime }\right) ^{\ast }\text{%
.}
\end{eqnarray*}%
We let $\Lambda ^{\prime }:z,z^{\prime }\in G\mapsto \Lambda (z^{\prime },z)$%
. Then an easy computation shows that $\Lambda ^{\prime }$ is a unitary
projective representation of $G$ on $\mathbb{C}^{2}$ associated to the
cocycle defined by $\exp \left( -i\pi \sigma \right) $, and $\Lambda
^{\prime }(\zeta )=V$ and $\Lambda ^{\prime }(\xi )=U$.

Let $B=A\rtimes _{\gamma }G$. Let us define the representation $\Psi $ of $B$
on $\mathbb{C}^{2}\otimes \mathbb{C}^{2}$ by:%
\begin{eqnarray*}
\Psi (a) &=&1\otimes a\text{,} \\
\Psi (U^{\zeta }) &=&V\otimes U\text{,} \\
\Psi (U^{\xi }) &=&U\otimes V\text{.}
\end{eqnarray*}%
First, we observe that:%
\begin{eqnarray*}
\Psi (U^{\zeta })\Psi (a)\Psi (U^{\zeta })^{\ast } &=&1\otimes UaU^{\ast
}=\Psi (\gamma _{\zeta }(a))\text{,} \\
\Psi (U^{\xi })\Psi (a)\Psi (U^{\xi })^{\ast } &=&1\otimes VaV^{\ast }=\Psi
(\gamma _{\xi }(a))\text{.}
\end{eqnarray*}%
Therefore $\Psi $ is indeed defining a representation of $B$. Moreover: 
\begin{equation*}
\Psi (U^{g})=\Lambda ^{\prime }(g)\otimes \Lambda (g)
\end{equation*}
for $g\in G$. Since $\Lambda ^{\prime }$ is irreducible, $\Psi $ is
irreducible as well by Theorem (\ref{Homogeneous}). Last, the commutant of $%
\Psi (A)$ is $M_{2}\left( \mathbb{C}\right) $, i.e. the restriction of $\Psi 
$ to $A$ is the direct sum of two copies of the identity representation of $A
$.
\end{example}

\bigskip We now turn to the special case of cyclic groups where the
representation $\Psi $ of Theorem (\ref{FiniteGroupConclusion}) is always
minimal, i.e. its restriction to $A$ is always an irreducible representation
of $A$. We shall characterize such minimal representations in terms of the
fixed point C*-subalgebra $A_{1}$ of $A$.

\section{Actions of Finite Cyclic Groups}

\bigskip Let $A$ be a unital C*-algebra and $\sigma $ be a *-automorphism of 
$A$ of period $n$, for $n\in \mathbb{N}$, i.e. $\sigma ^{n}=\limfunc{Id}_{A}$%
. We shall not assume that $n$ is the smallest such natural integer, i.e. $%
\sigma $ may be of an order dividing $n$. The automorphism $\sigma $
naturally generates an action of $\mathbb{Z}_{n}$ on $A$ by letting $\alpha
_{z}(a)=\sigma ^{k}(a)$ for all $z\in \mathbb{Z}_{n}$ and $k\in \mathbb{Z}$
of class $z$ modulo $n$. The crossed-product $A\rtimes _{\alpha }\mathbb{Z}%
_{n}$ will be simply denoted by $A\rtimes _{\sigma }\mathbb{Z}_{n}$, and the
canonical unitary $U^{1}\in A\rtimes _{\sigma }\mathbb{Z}_{n}$ corresponding
to $1\in \mathbb{Z}_{n}$ will simply be denoted by $U$. The C*-algebra $%
A\rtimes _{\sigma }\mathbb{Z}_{n}$ is universal among all C*-algebras
generated by a copy of $A$ and a unitary $u$ such that $u^{n}=1$ and $%
uau^{\ast }=\sigma (a)$.

\bigskip Theorem (\ref{FiniteGroupConclusion}) already provides much
information about the structure of irreducible representations of $A\rtimes
_{\sigma }\mathbb{Z}_{n}$. Yet we shall see it is possible in this case to
characterize these representations in terms of irreducible representations
of $A$ and of the fixed point C*-subalgebra $A_{1}$ of $A$ for $\sigma $. Of
central importance in this characterization are minimal representations of $A
$ for $\sigma $ and their relation to irreducible representations of $A_{1}$%
. We start this section with the exploration of this connection. Next, we
propose a full characterization of irreducible representations of $A\rtimes
_{\sigma }\mathbb{Z}_{n}$.

\subsection{Minimal Representations}

An extreme case of irreducible representation for crossed-products is given
by:

\begin{definition}
Let $\Pi $ be an irreducible representation of $A\rtimes _{\alpha }G$ is
called minimal when its restriction to $A$ is irreducible. Moreover, if $\pi 
$ is an irreducible representation of $A$ such that there exists some
irreducible representation $\Pi $ of $A\rtimes _{\alpha }G$ whose
restriction to $A$ is $\pi $, then we say that $\pi $ is minimal for the
action $\alpha $ of $G$.
\end{definition}

Such representations play a central role in the description of irreducible
representations of $A\rtimes _{\sigma }\mathbb{Z}_{n}$ when $\sigma $ is an
automorphism of period $n$. We propose to characterize them in term of the
fixed point C*-subalgebra $A_{1}$ of $A$. The set $\widehat{\mathbb{Z}_{n}}$
of irreducible representations of $\mathbb{Z}_{n}$ is the Pontryagin dual of 
$\mathbb{Z}_{n}$ which we naturally identify with the group $\mathbb{U}_{n}$
of $n^{\text{th}}$ roots of the unit in $\mathbb{C}$. Let $\lambda \in 
\mathbb{U}_{n}$. Thus $k\in \mathbb{Z}_{n}\mapsto \lambda ^{n}$ is an
irreducible representation of $\mathbb{Z}_{n}$ and the spectral subspace $%
A_{\lambda }$ of $A$ for $\lambda $ is given by $\left\{ a:\sigma
(a)=\lambda a\right\} $. Indeed, $A_{\lambda }$ is by definition the range
of the projection $P_{\lambda }:a\in A\mapsto \frac{1}{n}\sum_{k=0}^{n-1}%
\lambda ^{-k}\sigma ^{k}(a)$ by Equality (\ref{SpectralProjDef}), and it is
easy to check that $P_{\lambda }(a)=a\iff \sigma (a)=\lambda a$ from the
definition of $P_{\lambda }$.

\begin{theorem}
\bigskip \label{Rep}Let $\sigma $ be a *-automorphism of a unital C*-algebra 
$A$ of period $n$. Let $\Pi $ be an irreducible representation of $A\rtimes
_{\sigma }\mathbb{Z}_{n}$ on a Hilbert space $\mathcal{H}$ and let $\pi _{A}$
be its restriction to $A$. Let $\Sigma $ be the spectrum of $U_{\Pi }:=\pi
(U)$. Now, $\Sigma $ is a subset of $\mathbb{U}_{n}$; let us write $\Sigma
=\left\{ \lambda _{1},\ldots ,\lambda _{p}\right\} $ and denote the spectral
subspace of $U_{\Pi }$ associated to $\lambda _{j}$ by $\mathcal{H}_{j}$.
With the decomposition $\mathcal{H}=\oplus _{k=1}^{p}\mathcal{H}_{k}$, we
write, for all $a\in A$:%
\begin{equation}
\pi _{A}(a)=\left[ 
\begin{array}{cccc}
\alpha _{11}(a) & \alpha _{12}(a) & \cdots & \alpha _{1p}(a) \\ 
\alpha _{21}(a) & \alpha _{22}(a) &  & \alpha _{2p}(a) \\ 
\vdots &  & \ddots & \vdots \\ 
\alpha _{p1}(a) & \alpha _{p2}(a) & \cdots & \alpha _{nn}(a)%
\end{array}%
\right] \text{.}  \label{RepAlphaDec}
\end{equation}%
Then for $k,j\in \left\{ 1,\ldots ,p\right\} $ the map $\alpha _{jk}$ is a
linear map on $A_{\lambda _{j}\overline{\lambda _{k}}}$ and null on $\oplus
_{\mu \not=\lambda _{j}\overline{\lambda _{k}}}A_{\mu }$. Moreover, the maps 
$\alpha _{kk}$ are irreducible *-representations of the fixed point
C*-algebra $A_{1}$.

Furthermore, the following are equivalent:

\begin{itemize}
\item The representation $\pi _{A}$ of $A$ is irreducible, i.e. $\Pi $ is
minimal,

\item The *-representations $\alpha _{11},\ldots ,\alpha _{pp}$ are pairwise
not unitarily equivalent, i.e. for all $i\not=j\in \left\{ 1,\ldots
,p\right\} $ the representation $\alpha _{ii}$ is not equivalent to $\alpha
_{jj}$.
\end{itemize}
\end{theorem}

\begin{proof}
Since $U_{\Pi }^{n}=1$, the spectrum of the unitary $U_{\Pi }$ is a subset $%
\Sigma =\left\{ \lambda _{1},\ldots ,\lambda _{p}\right\} $ of $\mathbb{U}%
_{n}$ for some $p\in \mathbb{N}$. We write $\mathcal{H}=\mathcal{H}%
_{1}\oplus \ldots \oplus \mathcal{H}_{p}$ where $\mathcal{H}_{i}$ is the
spectral subspace of $U_{\Pi }$ for the eigenvalue $\lambda _{i}$ for $%
i=1,\ldots ,p$, so that $U_{\Pi }=\left[ 
\begin{array}{ccc}
\lambda _{1} &  &  \\ 
& \ddots &  \\ 
&  & \lambda _{p}%
\end{array}%
\right] $. Let $i,j\in \left\{ 1,\ldots ,p\right\} $ and let $\alpha _{ij}$
be the map defined by Identity (\ref{RepAlphaDec}). First, it is immediate
that $\alpha _{ij}$ is linear. Now, a simple computation shows that:%
\begin{eqnarray*}
&&U_{\Pi }\pi (a)U_{\Pi }^{\ast }= \\
&=&\left[ 
\begin{array}{ccc}
\lambda _{1} &  &  \\ 
& \ddots &  \\ 
&  & \lambda _{p}%
\end{array}%
\right] \left[ 
\begin{array}{ccc}
\alpha _{11}(a) & \cdots & \alpha _{1p}(a) \\ 
\vdots &  & \vdots \\ 
\alpha _{p1}(a) & \cdots & \alpha _{pp}(a)%
\end{array}%
\right] \left[ 
\begin{array}{ccc}
\overline{\lambda _{1}} &  &  \\ 
& \ddots &  \\ 
&  & \overline{\lambda _{p}}%
\end{array}%
\right] \\
&=&\left[ 
\begin{array}{cccc}
\alpha _{11}(a) & \lambda _{1}\overline{\lambda _{2}}\alpha _{12}(a) & \cdots
& \lambda _{1}\overline{\lambda _{p}}\alpha _{1p}(a) \\ 
\lambda _{2}\overline{\lambda _{1}}\alpha _{21}(a) & \alpha _{22}(a) &  & 
\lambda _{2}\overline{\lambda _{p}}\alpha _{2p}(a) \\ 
\vdots &  & \ddots & \vdots \\ 
\lambda _{p}\overline{\lambda _{1}}\alpha _{p1}(a) & \lambda _{p}\overline{%
\lambda _{2}}\alpha _{p2}(a) & \cdots & \alpha _{pp}(a)%
\end{array}%
\right] =\pi \left( \sigma (a)\right) \text{.}
\end{eqnarray*}

Therefore for all $i,j\in \left\{ 1,\ldots ,p\right\} $ we have that $\alpha
_{ij}(\sigma (a))=\lambda _{i}\overline{\lambda _{j}}\alpha _{ij}(a)$. Let $%
a\in A_{\mu }$ for $\mu \in \mathbb{U}_{n}$, i.e. $\sigma (a)=\mu a$. Then $%
\alpha _{ij}(\sigma (a))=\mu \alpha _{ij}(a)$. Therefore either $\alpha
_{ij}(a)=0$ or $\mu =\lambda _{i}\overline{\lambda _{j}}$.

In particular, $\alpha _{jj}$ is a representation of $A_{1}$ for all $j\in
\left\{ 1,\ldots ,p\right\} $. Indeed, if $a\in A_{1}$ then $\alpha
_{jk}(a)=0$ if $j\not=k$ and thus $\pi _{A}(a)$ is diagonal. Since $\pi _{A}$
is a representation of $A$, it follows from easy computations that $\alpha
_{jj}$ are representations of $A_{1}$.

Now, since $A\oplus AU\oplus \ldots \oplus AU^{n-1}=A\rtimes _{\sigma }%
\mathbb{Z}_{n}$, every element of the range of $\Pi $ is of the form $\oplus
_{j=0}^{n-1}\pi _{A}(a_{j})U_{\Pi }^{j}$ for $a_{0},\ldots ,a_{n-1}\in A$.
Now, let $i\in \left\{ 1,\ldots ,p\right\} $. We observe that the $(i,i)$
entry of $\oplus _{j=0}^{n-1}\pi _{A}(a_{i})U_{\Pi }^{j}$ in the
decomposition $\mathcal{H}=\mathcal{H}_{1}\oplus \ldots \oplus \mathcal{H}%
_{p}$ is given by $\sum_{j=0}^{n-1}\lambda _{i}^{j}\alpha
_{ii}(a_{j})=\alpha _{ii}\left( \sum_{j=0}^{n-1}\lambda _{i}^{j}a_{j}\right) 
$. Hence, the $(i,i)$ entries of operators in the range of $\Pi $ are
exactly given by the operators in the range of $\alpha _{ii}$. Now, let $T$
be any operator acting on $\mathcal{H}$. Since $\Pi $ is irreducible, by the
Von Neumann double commutant Theorem \cite[Theorem 1.7.1]{Davidson}, $T$ is
the limit, in the weak operator topology ($\limfunc{WOT}$), of elements in
the range of $\Pi $. In particular, the $(i,i)$ entry of $T$ in the
decomposition $\mathcal{H}=\mathcal{H}_{1}\oplus \ldots \oplus \mathcal{H}%
_{p}$ is itself a $\limfunc{WOT}$ limit of elements in the range of $\alpha
_{ii}$ since the left and right multiplications by a fixed operator are $%
\limfunc{WOT}$ continuous \cite[p. 16]{Davidson}. Therefore, the range of $%
\alpha _{ii}$ is $\limfunc{WOT}$ dense in $\mathcal{H}_{i}$. Thus, by the
double commutant theorem again, $\alpha _{ii}$ is irreducible.

We now turn to characterizing minimal representations. We first establish a
necessary condition.

Suppose that there exists $i,j\in \left\{ 1,\ldots ,p\right\} $ with $%
i\not=j $ and a unitary $u$ such that $u\alpha _{ii}u^{\ast }=\alpha _{jj}$.
In the decomposition $\mathcal{H}=\mathcal{H}_{1}\oplus \ldots \oplus 
\mathcal{H}_{p}$, define the block-diagonal unitary 
\begin{equation*}
D_{u}^{i}=\underset{i-1\text{ times}}{\underbrace{1\oplus \ldots \oplus 1}}%
\oplus u\oplus \underset{p-i\text{ times}}{\underbrace{1\oplus \ldots \oplus
1}}\text{.}
\end{equation*}

Then by conjugating $\pi _{A}$ by $D_{u}^{i}$, we see that we may as well
assume $\alpha _{ii}=\alpha _{jj}$. Yet, this implies that in the $\limfunc{%
WOT}$-closure of the range of $\pi _{A}$, every operator has the same $(i,i)$
and $(j,j)$ entry in the decomposition $\mathcal{H}=\mathcal{H}_{1}\oplus
\ldots \oplus \mathcal{H}_{p}$. Hence the range of $\pi _{A}$ is not $%
\limfunc{WOT}$-dense and thus $\pi _{A}$ is reducible, so $\Pi $ is not
minimal.

We now prove that our necessary condition is also sufficient. Assume that $%
\alpha _{11},\ldots ,\alpha _{pp}$ are pairwise not unitary equivalent. The
claim is that $\pi _{A}$ is irreducible.

Let $T\in \left( \pi \left( A\right) \right) ^{\prime }.$ Decompose $T=\left[
\begin{array}{ccc}
T_{11} & \cdots & T_{1p} \\ 
\vdots &  & \vdots \\ 
T_{p1} & \cdots & T_{pp}%
\end{array}%
\right] $ with respect to the decomposition $\mathcal{H}=\oplus _{i=1}^{p}%
\mathcal{H}_{i}$. Let $i\neq j.$ First, note that if $a\in A_{1}$ then $%
\alpha _{ij}(a)=0$. Second, since $T$ commutes with $\pi _{A}(a)$ for $a\in
A_{1}$, we have:%
\begin{equation}
\alpha _{ii}\left( a\right) T_{ij}=T_{ij}\alpha _{jj}\left( a\right) \text{
for }a\in A_{1}\text{.}  \label{Rep1}
\end{equation}

By Lemma (\ref{Schur}), since $\alpha _{ii}$ and $\alpha _{jj}$ are
irreducible and not unitarily equivalent for $i\not=j$, we conclude that $%
T_{ij}=0$. Moreover, for all $i\in \left\{ 1,\ldots ,p\right\} $ and $a\in
A_{1}$ we have $\alpha _{ii}\left( a\right) T_{ii}=T_{ii}\alpha _{ii}\left(
a\right) $. Since $\alpha _{ii}$ is irreducible, we conclude that $T_{ii}$
is a scalar. Therefore, the operator $T$ commutes with the operator $U_{\Pi
} $. Since $\Pi $ is irreducible, we conclude that $T$ itself is a scalar.
Therefore, $\pi _{A}$ is an irreducible representation of $A$ and thus $\Pi $
is minimal.
\end{proof}

\bigskip Together with Theorem (\ref{FiniteGroupConclusion}), Theorem (\ref%
{Rep}) will allow us to now develop further the description of arbitrary
irreducible representations of crossed-products by finite cyclic groups. It
is interesting to look at a few very simple examples to get some intuition
as to what could be a more complete structure theory for irreducible
representations of crossed-products by $\mathbb{Z}_{n}$. First of all, one
should not expect in general that the spectrum of $U_{\Pi }$ is a coset of $%
\mathbb{Z}_{n}$, as the simple action of $\sigma =\limfunc{Ad}\left[ 
\begin{array}{cc}
i &  \\ 
& e^{i\frac{3\pi }{4}}%
\end{array}%
\right] $ on $M_{2}\left( \mathbb{C}\right) $ shows. In this case, the
identity is the only irreducible representation of the crossed-product $%
M_{2}\left( \mathbb{C}\right) \rtimes _{\sigma }\mathbb{Z}_{4}=M_{2}\left( 
\mathbb{C}\right) $ and clearly $\left\{ i,e^{i\frac{3\pi }{4}}\right\} $ is
not a coset of $\mathbb{Z}_{4}$. Of course, this is an example of a minimal
representation.

\bigskip In \cite{Latremoliere06}, we showed that all irreducible
representations of $A\rtimes _{\sigma }\mathbb{Z}_{2}$ where regular or
minimal. The following example shows that we can not expect the same in the
general case.

\begin{example}
\label{CuteExample}Let $A=M_{2}\left( \mathbb{C}\right) \oplus M_{2}\left( 
\mathbb{C}\right) $ and define $\sigma (M\oplus N)=WNW^{\ast }\oplus M$ with 
$W=\left[ 
\begin{array}{cc}
0 & 1 \\ 
1 & 0%
\end{array}%
\right] $. Then $\sigma ^{4}=\limfunc{Id}_{A}$ and $\sigma ^{2}(M\oplus
N)=WMW^{\ast }\oplus WNW^{\ast }$. Now, let $\pi _{i}:M_{1}\oplus M_{2}\in
A\mapsto M_{i}$ with $i=1,2$. Of course, $\pi _{1},\pi _{2}$ are the only
two irreducible representations of $A$ up to equivalence, and they are not
equivalent to each other (since they have complementary kernels). Now, we
consider the following representation $\Pi $ of $A\rtimes _{\sigma }\mathbb{Z%
}_{4}$. It acts on $\mathbb{C}^{4}$. We set:%
\begin{equation*}
\pi _{A}=\left[ 
\begin{array}{cc}
\pi _{1} & 0 \\ 
0 & \pi _{2}%
\end{array}%
\right]
\end{equation*}%
and:%
\begin{equation*}
U_{\Pi }=\left[ 
\begin{array}{cc}
0 & 1 \\ 
W & 0%
\end{array}%
\right] \text{.}
\end{equation*}%
First, observe that $\Pi $ thus defined is irreducible. Indeed, $M$ commutes
with $\pi _{A}$ if and only if $M=\left[ 
\begin{array}{cc}
\lambda & b \\ 
c & \mu%
\end{array}%
\right] $ with $\lambda ,\mu \in \mathbb{C}$ and $b\pi _{2}(a)=\pi _{1}(a)c$
with $a\in A$. Now, $M$ commutes with $U_{\Pi }$ if and only if $\lambda
=\mu $ and $Wb=c$. Now, let $a\in M_{2}\left( \mathbb{C}\right) $ be
arbitrary; then $b\pi _{2}\left( a\oplus Wa\right) =\pi _{1}\left( a\oplus
Wa\right) c$ i.e.%
\begin{equation*}
bWa=abW\text{.}
\end{equation*}%
Hence $bW$ is scalar. So $b=\lambda W$. Thus $b$ commutes with $W$. But then
for an arbitrary $a$ we have $b\pi _{2}\left( aW\oplus a\right) =\pi
_{1}\left( aW\oplus a\right) bW$ i.e. $ba=aWbW=ab$ so $b$ commutes with $%
M_{2}\left( \mathbb{C}\right) $ and thus is scalar. Hence $b=0$. So $%
M=\lambda 1$ for $\lambda \in \mathbb{C}$ as needed.

Moreover, the restriction of $\Pi $ to $A$ is $\pi _{A}=\pi _{1}\oplus \pi
_{2}$. Thus, $\pi _{A}$ is reducible. Now, the fixed point C*-algebra $A_{1}$
is the C*-algebra $\left\{ M\oplus M:M=\left[ 
\begin{array}{cc}
a & b \\ 
b & a%
\end{array}%
\right] ;a,b\in \mathbb{C}\right\} $. Thus, $A_{1}$ has two irreducible
representations which are not equivalent:%
\begin{equation*}
\varphi _{1}:\left[ 
\begin{array}{cc}
a & b \\ 
b & a%
\end{array}%
\right] \in A_{1}\mapsto a+b
\end{equation*}%
and%
\begin{equation*}
\varphi _{2}:\left[ 
\begin{array}{cc}
a & b \\ 
b & a%
\end{array}%
\right] \in A_{1}\mapsto a-b\text{.}
\end{equation*}

We note that for $i=1,2$ we have $\pi _{i}$ restricted to $A_{1}$ is $%
\varphi _{1}\oplus \varphi _{2}$.
\end{example}

Now, using the notations of Example (\ref{CuteExample}), $\Pi $ is not
regular, since the restriction of any irreducible regular representation to
the fixed point algebra $A_{1}$ is given by the sum of several copies of the
same irreducible representation of $A_{1}$. Trivially, $\Pi $ is not minimal
either since $\Pi _{|A}=\pi _{1}\oplus \pi _{2}$. However, both $\pi _{1}$
and $\pi _{2}$ are minimal for the action of $\sigma ^{2}$. Moreover, both $%
\pi _{1}$ and $\pi _{2}$ restricted to $A_{1}$ are the same representation $%
\alpha _{1}\oplus \alpha _{2}$. We shall see in the next section that this
pattern is in fact general.

\subsection{Characterization of Irreducible Representations}

We now present the main result of this paper concerning crossed products by
finite cyclic groups. In this context, one can go further than Theorem\ (\ref%
{FiniteGroupConclusion}) to obtain a characterization of irreducible
representations of the crossed-products in term of the C*-algebras $A$ and $%
A_{1}$. The next lemma is the sufficient condition for this characterization.

\begin{lemma}
\label{SufficientCyclic}Let $\pi _{1}$ be an irreducible representation of $%
A $ acting on a Hilbert space $\mathcal{J}$. Assume that there exists a
unitary $V$ on $\mathcal{J}$ such that for some $m,k\in \left\{ 1,\ldots
,n\right\} $ with $n=mk$ we have $\pi _{1}\circ \sigma ^{m}=V\pi _{1}V^{\ast
}$ and $V^{k}=1$, and that $m$ is the smallest such nonzero natural integer,
i.e. $\pi _{1}\circ \sigma ^{j}$ is not unitarily equivalent to $\pi _{1}$
for $j\in \left\{ 2,\ldots ,m-1\right\} $. Then define the following
operators on the Hilbert space $\mathcal{H}=\underset{m\text{ times}}{%
\underbrace{\mathcal{J}\oplus \ldots \oplus \mathcal{J}}\text{:}}$%
\begin{equation*}
\Pi \left( U\right) =\left[ 
\begin{array}{ccccc}
0 & 1 & 0 & \cdots & 0 \\ 
\vdots &  & 1 &  & \vdots \\ 
\vdots &  &  & \ddots & \vdots \\ 
0 & 0 & \cdots & 0 & 1 \\ 
V & 0 & \cdots & 0 & 0%
\end{array}%
\right]
\end{equation*}%
and for all $a\in A$:%
\begin{equation*}
\pi _{A}(a)=\left[ 
\begin{array}{cccc}
\pi _{1}(a) &  &  &  \\ 
& \pi _{1}\circ \sigma (a) &  &  \\ 
&  & \ddots &  \\ 
&  &  & \pi _{1}\circ \sigma ^{m-1}(a)%
\end{array}%
\right] \text{.}
\end{equation*}%
Then the unique extension of $\Pi $ to $A\rtimes _{\sigma }\mathbb{Z}_{n}$
is an irreducible representation of $A\rtimes _{\sigma }\mathbb{Z}_{n}$.
\end{lemma}

\begin{proof}
An easy computation shows that $\Pi $ thus defined is a representation of $%
A\rtimes _{\sigma }\mathbb{Z}_{n}$ on $\mathcal{H}=\underset{m\text{ times}}{%
\underbrace{\mathcal{J}\oplus \ldots \oplus \mathcal{J}}}$. Write $\pi
_{i}=\pi _{1}\circ \sigma ^{i-1}$ for $i=1,\ldots ,m$. Let $T$ be an
operator which commutes with the range of $\Pi $. Then $T$ commutes with $%
\pi _{A}:=\Pi _{|A}$. Writing $T$ in the decomposition $\mathcal{H}=\mathcal{%
J}\oplus \ldots \oplus \mathcal{J}$ as:%
\begin{equation*}
T=\left[ 
\begin{array}{ccc}
T_{11} & \cdots & T_{1m} \\ 
\vdots &  & \vdots \\ 
T_{m1} & \cdots & T_{mm}%
\end{array}%
\right]
\end{equation*}%
Let $i,j\in \left\{ 1,\ldots ,m\right\} $. Since $T\pi _{A}(a)=\pi _{A}(a)T$
for all $a\in A$, we conclude that $\pi _{i}(a)T_{ij}=T_{ij}\pi _{j}(a)$. By
Lemma\ (\ref{Schur}), since $\pi _{i}$ and $\pi _{j}$ are irreducible and
not unitarily equivalent, we conclude that $T_{ij}=0$. Moreover, $T_{ii}$
commutes with $\pi _{i}$ which is irreducible, so we conclude that:%
\begin{equation*}
T=\left[ 
\begin{array}{ccc}
\lambda _{1} &  &  \\ 
& \ddots &  \\ 
&  & \lambda _{m}%
\end{array}%
\right]
\end{equation*}%
for $\lambda _{1},\ldots ,\lambda _{m}\in \mathbb{C}$. Since $T$ commutes
with $U_{\Pi }$ we conclude that $\lambda _{1}=\lambda _{i}$ for all $i\in
\left\{ 1,\ldots ,m\right\} $. Hence $\Pi $ is irreducible.
\end{proof}

\bigskip We now are ready to describe in detail the structure of irreducible
representations of crossed-products by finite cyclic groups in terms of
irreducible representations of $A$ and $A_{1}$.

\begin{theorem}
\label{CyclicConclusion}Let $\sigma $ be a *-automorphism of period $n$ of a
unital C*-algebra $A$. Then the following are equivalent:

\begin{enumerate}
\item $\Pi $ is an irreducible representation of $A\rtimes _{\sigma }\mathbb{%
Z}_{n}$,

\item There exists $k,m\in \mathbb{N}$ with $km=n$, an irreducible
representation $\pi _{1}$ of $A$ on a Hilbert space $\mathcal{J}$ and a
unitary $V$ on $\mathcal{J}$ such that $V^{k}=1$ and $V\pi _{1}\left( \cdot
\right) V=\pi _{1}\circ \sigma ^{m}\left( \cdot \right) $ such that:%
\begin{equation*}
\Pi (U)=\left[ 
\begin{array}{cccc}
0 & 1 &  &  \\ 
& \ddots & \ddots &  \\ 
&  & 0 & 1 \\ 
V &  &  & 0%
\end{array}%
\right]
\end{equation*}%
and for all $a\in A$:%
\begin{equation*}
\Pi (a)=\left[ 
\begin{array}{cccc}
\pi _{1}(a) &  &  &  \\ 
& \pi _{1}\circ \sigma (a) &  &  \\ 
&  & \ddots &  \\ 
&  &  & \pi _{1}\circ \sigma ^{m-1}(a)%
\end{array}%
\right]
\end{equation*}%
where for any $i\in \left\{ 1,\ldots ,m-1\right\} $ the representations $\pi
_{1}$ and $\pi _{1}\circ \sigma ^{i}$ are not equivalent.
\end{enumerate}

Moreover, if (2) holds then the representation $\psi $ of $A\rtimes _{\sigma
^{m}}\mathbb{Z}_{k}$ on $\mathcal{J}$\ defined by $\psi (a)=\pi _{1}(a)$ for 
$a\in A$ and $\psi (U)=V$ is a minimal representation of $A\rtimes _{\sigma
^{m}}\mathbb{Z}_{k}$. Let $\eta $ be the cardinal of the spectrum of $V$.
The restriction of $\pi _{1}$ to $A_{1}$ is therefore the sum of $\eta $
irreducible representations $\varphi _{1},\ldots ,\varphi _{\eta }$ of $%
A_{1} $ which are not pairwise equivalent. Last, the restriction of $\pi
_{1}\circ \sigma ^{i}$ to $A_{1}$ is unitarily equivalent to $\varphi
_{1}\oplus \ldots \oplus \varphi _{\eta }=\pi _{1|A_{1}}$ for all $i\in
\left\{ 0,\ldots ,m-1\right\} $.
\end{theorem}

\begin{proof}
By Lemma\ (\ref{SufficientCyclic}), (2) implies (1). We now turn to the
proof of (1) implies (2). Let $\Pi $ be an irreducible representation of $%
A\rtimes _{\sigma }\mathbb{Z}_{n}$. By Theorem (\ref{FiniteGroupConclusion}%
), there exists $m\in \mathbb{N}$ such that $m$ divides $n$, an irreducible
representation $\pi _{1}$ of $A$ on some space $\mathcal{H}_{1}$ and $r\in 
\mathbb{N}$ with $r>0$ such that, if $\pi =r\cdot \pi _{1}$ then up to
conjugating $\Pi $ by some unitary:

\begin{itemize}
\item For all $i=1,\ldots ,m-1$ the representation $\pi \circ \sigma ^{i}$
is not equivalent to $\pi $,

\item The representation $\pi \circ \sigma ^{m}$ is equivalent to $\pi $,

\item We have the decomposition $\mathcal{H}=\mathcal{J}_{0}\oplus \ldots
\oplus \mathcal{J}_{m-1}$ where $\mathcal{J}_{i}$ is the space on which $%
\left( r\cdot \pi \right) \circ \sigma ^{i}$ acts for $i\in \left\{ 0,\ldots
,m\right\} $ and is isometrically isomorphic to $\mathcal{J}$,

\item In the decomposition, $\mathcal{H}=\mathcal{J}_{0}\oplus \ldots \oplus 
\mathcal{J}_{m-1}$ there exists unitaries $U_{1},\ldots ,U_{m}$ such that:%
\begin{equation*}
U_{\Pi }=\left[ 
\begin{array}{ccccc}
0 & U_{1} & 0 & \cdots & 0 \\ 
0 & 0 & U_{2} & 0 & \vdots \\ 
\vdots &  & \ddots & \ddots & 0 \\ 
0 &  &  & 0 & U_{m-1} \\ 
U_{m} & 0 & \cdots & 0 & 0%
\end{array}%
\right]
\end{equation*}%
with $\left( \pi _{i}\circ \sigma \right) U_{i}=U_{i}\pi _{i+1}$ and $U_{i}:%
\mathcal{H}_{i+1}\longrightarrow \mathcal{H}_{i}$ for all $i\in \mathbb{Z}%
_{m}$.
\end{itemize}

Indeed, if $G=\mathbb{Z}_{n}$ in Theorem (\ref{FiniteGroupConclusion}) then $%
H$, as a subgroup of $G$, is of the form $\left( m\mathbb{Z}\right) /n%
\mathbb{Z}$ with $m$ dividing $n$, and if we let $g_{1}=0$, $g_{2}=1$,
\ldots , $g_{m}=m-1$ then we can check that this choice satisfies the
hypothesis of Theorem (\ref{FiniteGroupConclusion}). With this choice, the
permutation $\sigma ^{1}$ is then easily seen to be given by the cycle $%
\left( 1~2~\ldots ~m\right) $.

We will find it convenient to introduce some notation for the rest of the
proof. By Theorem (\ref{FiniteGroupConclusion}), for $i\in \left\{ 0,\ldots
,m-1\right\} $, after possibly conjugating $\Pi $ by some unitary, we can
decompose $\mathcal{J}_{i}$ as $\mathcal{H}_{ri+1}\oplus \ldots \oplus 
\mathcal{H}_{r(i+1)}$, where $\mathcal{H}_{ri+j}$ is isometrically
isomorphic to $\mathcal{H}_{1}$ for all $j\in \left\{ 1,\ldots ,r\right\} $,
so that the restriction of $\Pi _{|A}$ to the space $\mathcal{J}_{i}$ is
written $\left( \underset{r\text{ times}}{\underbrace{\pi _{1}\oplus \ldots
\oplus \pi _{1}}}\right) \circ \sigma ^{i}$ in this decomposition.

We now show how to conjugate $\Pi $ by a unitary to simplify its expression
further.

If we define the unitary $\Upsilon $ from $\mathcal{H}=\mathcal{J}_{0}\oplus
\ldots \oplus \mathcal{J}_{m-1}$ onto $\oplus _{1}^{m}\mathcal{J}_{m-1}$ by: 
\begin{equation*}
\Upsilon =\left[ 
\begin{array}{cccc}
U_{m}^{\ast }U_{m-1}^{\ast }\cdots U_{1}^{\ast } &  &  &  \\ 
& U_{m}^{\ast }\cdots U_{2}^{\ast } &  &  \\ 
&  & \ddots &  \\ 
&  &  & U_{m}^{\ast }%
\end{array}%
\right]
\end{equation*}%
then the unitary $\limfunc{Ad}\left( \Upsilon \right) \circ \Pi \left(
U\right) $ of $\oplus _{1}^{m}\mathcal{J}_{m-1}$\ is of the simpler form%
\begin{equation}
\limfunc{Ad}\Upsilon \circ \Pi \left( U\right) =\left[ 
\begin{array}{ccccc}
0 & 1 & 0 & \cdots & 0 \\ 
0 & 0 & 1 & \ddots & \vdots \\ 
\vdots & \vdots & \ddots & \ddots & 0 \\ 
0 & 0 & \cdots & 0 & 1 \\ 
V & 0 & \cdots & 0 & 0%
\end{array}%
\right]  \label{UnitaryShift}
\end{equation}%
for some unitary $V$ of $\mathcal{J}_{m-1}$. Moreover, if we write $\rho
_{1}=\limfunc{Ad}\left( U_{i}^{\ast }\ldots U_{1}^{\ast }\right) \circ \pi
_{1}$, then:%
\begin{equation*}
\limfunc{Ad}\Upsilon \circ \pi _{A}=\dbigoplus\limits_{j=1}^{m}\left( 
\underset{r\text{ times}}{\underbrace{\rho _{1}\circ \sigma ^{j-1}\oplus
\ldots \oplus \rho _{1}\circ \sigma ^{j-1}}}\right)
\end{equation*}%
and $\rho _{1}$ is by definition an irreducible representation of $A$
unitarily equivalent to $\pi _{1}$.

To simplify notations, we shall henceforth drop the notation $\limfunc{Ad}%
\Upsilon $ and simply write $\Pi $ for $\limfunc{Ad}\Upsilon \circ \Pi $. In
other words, we replace $\Pi $ by $\limfunc{Ad}\Upsilon \circ \Pi $ and we
shall use the notations introduced to study $\Pi $ henceforth, with the
understanding that for all $j=0,\ldots ,m-1$ and $k=1,\ldots ,r$ we have
that $\pi _{rj+k}=\pi _{1}\circ \sigma ^{j}$, that $\mathcal{J}_{j}$ is an
isometric copy of $\mathcal{J}_{0}$ and that $\mathcal{H}$ $=\mathcal{J}%
_{0}\oplus \ldots \oplus \mathcal{J}_{m-1}$ with $\mathcal{J}_{j}=\mathcal{H}%
_{rj+1}\oplus \ldots \oplus \mathcal{H}_{r(j+1)}$ where $\pi _{rj+k}$ acts
on $\mathcal{H}_{rj+k}$ which is an isometric copy of $\mathcal{H}_{1}$.
Moreover, $U_{\Pi }$ is given by Equality (\ref{UnitaryShift})\ for some
unitary $V$ of $\mathcal{J}_{0}$.

We are left to show that each irreducible subrepresentation of $\pi _{A}$ is
of multiplicity one, i.e. $r=1$. We recall that we have shown above that $%
H=\left( m\mathbb{Z}\right) /n\mathbb{Z}$ with $n=mk$ and $k\in \mathbb{N}$.
Using the notations of Theorem (\ref{FiniteGroupConclusion}), the
representation $\Psi $ defined by $\Psi (a)=\pi (a)$ for all $a\in A$ and $%
\Psi (U^{m})=V$ is an irreducible representation of $A\rtimes _{\alpha }H$.
Now $A\rtimes _{\alpha }H$ is *-isomorphic to $A\rtimes _{\alpha ^{m}}%
\mathbb{Z}_{k}$ by universality of the C*-crossed-product, and we now
identify these two C*-algebras. The image of $U^{m}\in A\rtimes _{\alpha }H$
in the crossed-product $A\rtimes _{\alpha ^{m}}\mathbb{Z}_{k}$ is denoted by 
$\upsilon $ and is the canonical unitary of $A\rtimes _{\alpha ^{m}}\mathbb{Z%
}_{k}$. Thus by Theorem (\ref{FiniteGroupConclusion}) $\Psi $ is an
irreducible representation of $A\rtimes _{\alpha ^{m}}\mathbb{Z}_{k}$ which
(up to conjugacy) acts on the space $\mathbb{C}^{r}\otimes \mathcal{H}_{1}$
and is of the form $\Psi (a)=1_{\mathbb{C}^{r}}\otimes \pi _{1}(a)$ for $%
a\in A$ and $\Psi (\upsilon ^{z})=\Omega (z)\otimes W(z)$ for $z\in \mathbb{Z%
}_{k}$ where $\Omega $ and $W$ are some unitary projective representations
of $\mathbb{Z}_{k}$ on $\mathbb{C}^{r}$ and $\mathcal{H}_{1}$ respectively,
with $\Omega $ being irreducible. Since $\mathbb{Z}_{k}$ is cyclic, the
range of the projective representation $\Omega $ is contained in the
C*-algebra $C^{\ast }\left( \Omega (1)\right) $ which is Abelian since $%
\Omega (1)$ is a unitary. Hence, since $\Omega $ is irreducible, $C^{\ast
}\left( \Omega (1)\right) $ is an irreducible Abelian C*-algebra of
operators acting on $\mathbb{C}^{r}$. Hence $r=1$ and $\mathcal{J}=\mathcal{H%
}_{1}$. Moreover, since $U_{\Pi }^{m}=V\oplus \ldots \oplus V$ then $U_{\Pi
}^{n}=V^{k}\oplus \ldots \oplus V^{k}=1_{\mathcal{H}}$ and thus $V^{k}=1_{%
\mathcal{J}}$. Therefore, (2) holds as claimed.

Last, we also observed that $V\pi _{1}V=\pi _{1}\circ \sigma ^{k}$ by
construction (since $U_{\Pi }^{k}=V\oplus \ldots \oplus V$). Hence by
definition, since $\pi _{1}$ is irreducible, the representation $\psi $ of $%
A\rtimes _{\sigma ^{k}}\mathbb{Z}_{\mu }$ defined by $\psi (a)=\pi _{1}(a)$
for $a\in A$ and $\psi (U)=V$ is minimal. Hence, by Theorem (\ref{Rep}), the
restriction of $\pi _{1}$ to the fixed point C*-algebra $A_{1}$ is the
direct sum of $\eta $ irreducible representations $\varphi _{1},\ldots
,\varphi _{\eta }$ of $A_{1}$ such that $\varphi _{i}$ and $\varphi _{j}$
are not unitarily equivalent for $i\not=j\in \left\{ 1,\ldots ,\eta \right\} 
$, where $\eta $ is the cardinal of the spectrum of $V$. Moreover, since $%
\pi _{i}=\pi _{1}\circ \sigma ^{i}$ it is immediate that $\pi _{i}$
restricted to $A_{1}$ equals to $\pi _{1}$ restricted to $A_{1}$. This
concludes our proof.
\end{proof}

\begin{corollary}
Let $\Pi $ be an irreducible representation of $A\rtimes _{\sigma }\mathbb{Z}%
_{n}$. The following are equivalent:

\begin{enumerate}
\item Up to unitary equivalence, $\Pi $ is an irreducible regular
representation of $A\rtimes _{\sigma }\mathbb{Z}_{n}$, i.e. it is induced by
a unique irreducible representation $\pi $ of $A$ and:%
\begin{equation*}
U_{\Pi }=\left[ 
\begin{array}{cccc}
0 & 1 &  &  \\ 
& 0 & \ddots &  \\ 
&  & \ddots & 1 \\ 
1 &  &  & 0%
\end{array}%
\right]
\end{equation*}%
while $\pi _{A}=\oplus _{i=0}^{n-1}\pi \circ \sigma ^{i}$ and $\pi \circ
\sigma ^{i}$ is not equivalent to $\pi \circ \sigma ^{j}$ for $i,j=1,\ldots
,n-1$ with $i\not=j$,

\item There exists an irreducible subrepresentation $\pi $ of $\Pi _{|A}$
such that $\pi \circ \sigma ^{i}$ is not equivalent to $\pi $ for $%
i=1,\ldots ,n-1$,

\item There exists a unique irreducible representation $\varphi $ of $A_{1}$
such that $\Pi _{|A_{1}}$ is equivalent to $n\cdot \varphi $,

\item There is no $k\in \left\{ 1,\ldots ,n-1\right\} $ such that the
C*-algebra generated by $\Pi (A)$ and $U_{\Pi }^{k}$ is reducible.
\end{enumerate}
\end{corollary}

\begin{proof}
It is a direct application of Theorem\ (\ref{CyclicConclusion}).
\end{proof}

We thus have concluded that all irreducible representations of crossed
products by finite cyclic groups have a structure which is a composite of
the two cases found in \cite{Latremoliere06}. Indeed, such representations
cycle through a collection of minimal representations, which all share the
same restriction to the fixed point algebra. The later is a finite sum of
irreducible mutually disjoint representations of the fixed point algebra.

\begin{remark}
Let $\sigma $ be an order $n$ automorphism of a unital C*-algebra $A$ and
let $\Pi $ be an irreducible representation of $A\rtimes _{\sigma }\mathbb{Z}
$. We recall \cite{Zeller-Meier68} that $A\rtimes _{\sigma }\mathbb{Z}$ is
generated by $A$ and a unitary $U$ such that $UaU^{\ast }=\sigma (a)$ for
all $a\in A$ and is universal for these commutation relations. We denote $%
\Pi (U)$ by $U_{\Pi }$ and $\Pi (a)$ by $\pi (a)$ for all $a\in A$. Now,
note that $U_{\Pi }^{n}$ commutes with $\pi $ since $\sigma ^{n}=\limfunc{Id}%
_{A}$ and of course $U_{\Pi }^{n}$ commutes with $U_{\Pi }$ so, since $\Pi $
is irreducible, there exists $\lambda \in \mathbb{T}$ such that $U_{\Pi
}^{n}=\lambda $. Now, define $V_{\Pi }=\overline{\mu }U_{\Pi }$ for any $\mu
\in \mathbb{T}$ such that $\mu ^{n}=\lambda $. Then $V_{\Pi }^{n}=1$ and
thus $\left( \pi ,V_{\Pi }\right) $ is an irreducible representation of $%
A\rtimes _{\sigma }\mathbb{Z}_{n}$ which is then fully described by Theorem (%
\ref{CyclicConclusion}).
\end{remark}

In the last section of this paper, we give a necessary condition on
irreducible representations of crossed-products by the group $\mathfrak{S}%
_{3}$ of permutations of $\left\{ 1,2,3\right\} $. This last example
illustrates some of the behavior which distinguish the conclusion of Theorem
(\ref{FiniteGroupConclusion}) from the one of Theorem (\ref{CyclicConclusion}%
).

\section{Application: Crossed-Products by the permutation group on $\left\{
1,2,3\right\} $}

As an application, we derive the structure of the irreducible
representations of crossed-products by the group $\mathfrak{S}_{3}$ of
permutations of $\left\{ 1,2,3\right\} $. This group is isomorphic to $%
\mathbb{Z}_{3}\rtimes _{\gamma }\mathbb{Z}_{2}$ where $\gamma $ is defined
as follows:\ if $\eta $ and $\tau $ are the respective images of $1\in 
\mathbb{Z}$ in the groups $\mathbb{Z}_{3}$ and $\mathbb{Z}_{2}$ then the
action $\gamma $ of $\mathbb{Z}_{2}$ on $\mathbb{Z}_{3}$ is given by $\gamma
_{\tau }(\eta )=\eta ^{2}$. Thus in $\mathbb{Z}_{3}\rtimes _{\gamma }\mathbb{%
Z}_{2}$ we have $\tau \eta \tau =\eta ^{2}$, $\tau ^{2}=1$ and $\eta ^{3}=1$
(using the multiplicative notation for the group law). An isomorphism
between $\mathfrak{S}_{3}$ and $\mathbb{Z}_{3}\rtimes _{\gamma }\mathbb{Z}%
_{2}$ is given by sending the transposition $\left( 1~2\right) $ to $\tau $
and the $3$-cycle $\left( 1~2~3\right) $ to $\eta $. From now on we shall
identify these two groups implicitly using this isomorphism.

\begin{theorem}
\label{Permutation3}Let $\alpha $ be an action of $\mathfrak{S}_{3}$ on $A$.
Let $\Pi $ be an irreducible representation of $A\rtimes _{\alpha }\mathfrak{%
S}_{3}$. We denote by $\tau $ and $\eta $ the permutations $\left(
1~2\right) $ and $\left( 1~2~3\right) $. The set $\left\{ \tau ,\eta
\right\} $ is a generator set of $\mathfrak{S}_{3}$. We denote by $U_{\tau }$
and $U_{\eta }$ the canonical unitaries in $A\rtimes _{\alpha }\mathfrak{S}%
_{3}$ corresponding respectively to $\tau $ and $\eta $. Then either (up to
a unitary conjugation of $\Pi $):

\begin{itemize}
\item $\Pi $ is minimal, i.e. $\Pi _{|A}$ is irreducible,

\item There exists an irreducible representation $\pi _{1}$ on $\mathcal{H}%
_{1}$ of $A$ such that $\mathcal{H}=\mathcal{H}_{1}\oplus \mathcal{H}_{1}$
with $\pi _{A}=\pi _{1}\oplus \pi _{1}\circ \alpha _{\tau }$. Then $\Pi
(U_{\tau })=\left[ 
\begin{array}{cc}
0 & 1 \\ 
1 & 0%
\end{array}%
\right] $ in this decomposition. \emph{Observe that }$\pi _{1}$ \emph{may or
not be equivalent to} $\pi _{1}\circ \alpha _{\tau }$. Moreover, $\pi _{1}$
and $\pi _{1}\circ \alpha _{\tau }$ are minimal for the action of $\eta $.

\item There exists an irreducible representation $\pi _{1}$ on $\mathcal{H}%
_{1}$ of $A$ such that $\pi _{1}$ and $\pi _{1}\circ \alpha _{\eta ^{i}}$
are non equivalent for $i=1,2$ and such that $\mathcal{H}=\mathcal{H}%
_{1}\oplus \mathcal{H}_{1}\oplus \mathcal{H}_{1}$ with $\pi _{A}=\pi
_{1}\oplus \pi _{1}\circ \alpha _{\eta }\oplus \pi _{1}\circ \alpha _{\eta
^{2}}$. Then $\Pi (U_{\eta })=\left[ 
\begin{array}{ccc}
0 & 1 & 0 \\ 
0 & 0 & 1 \\ 
1 & 0 & 0%
\end{array}%
\right] $ in this decomposition.

\item Last, there exists an irreducible representation $\pi _{1}$ on $%
\mathcal{H}_{1}$ of $A$ such that $\pi _{1}\circ \alpha _{\sigma }$ is not
equivalent to $\pi _{1}$ for $\sigma \in \mathfrak{S}_{3}\backslash \left\{ 
\limfunc{Id}\right\} $ and $\mathcal{H}=\mathcal{H}_{1}^{\oplus 6}$ with:%
\begin{equation*}
\pi _{A}=\pi _{1}\oplus \pi _{1}\circ \alpha _{\eta }\oplus \pi _{1}\circ
\alpha _{\eta ^{2}}\oplus \pi _{1}\circ \alpha _{\tau }\oplus \pi _{1}\circ
\alpha _{\eta \tau }\oplus \pi _{1}\circ \alpha _{\eta ^{2}\tau }
\end{equation*}%
and%
\begin{equation*}
\Pi \left( U_{\eta }\right) =\left[ 
\begin{array}{cccccc}
0 & 1 & 0 & 0 & 0 & 0 \\ 
0 & 0 & 1 & 0 & 0 & 0 \\ 
1 & 0 & 0 & 0 & 0 & 0 \\ 
0 & 0 & 0 & 0 & 0 & 1 \\ 
0 & 0 & 0 & 1 & 0 & 0 \\ 
0 & 0 & 0 & 0 & 1 & 0%
\end{array}%
\right] \text{,}
\end{equation*}%
while%
\begin{equation*}
\Pi \left( U_{\tau }\right) =\left[ 
\begin{array}{cccccc}
0 & 0 & 0 & 1 & 0 & 0 \\ 
0 & 0 & 0 & 0 & 1 & 0 \\ 
0 & 0 & 0 & 0 & 0 & 1 \\ 
1 & 0 & 0 & 0 & 0 & 0 \\ 
0 & 1 & 0 & 0 & 0 & 0 \\ 
0 & 0 & 1 & 0 & 0 & 0%
\end{array}%
\right] \text{.}
\end{equation*}
\end{itemize}
\end{theorem}

\begin{proof}
The C*-algebra $A\rtimes _{\alpha }\mathfrak{S}_{3}$ is generated by a copy
of $A$ and two unitaries $U_{\tau }$ and $U_{\eta }$ that satisfy $U_{\tau
}^{2}=U_{\eta }^{3}=1$, $U_{\tau }U_{\eta }U_{\tau }=U_{\eta }^{2}$ and for
all $a\in A$ we have $U_{\tau }aU_{\tau }^{\ast }=\alpha _{\tau }(a)$ and $%
U_{\eta }aU_{\eta }^{\ast }=\alpha _{\eta }(a)$. Notice that $\mathfrak{S}%
_{3}=\mathbb{Z}_{3}\rtimes _{\gamma }\mathbb{Z}_{2}$ with $\gamma _{\tau
}\left( \eta \right) =\tau \eta \tau $. So we have $A\rtimes _{\alpha }%
\mathfrak{S}_{3}=\left( A\rtimes _{\alpha _{\eta }}\mathbb{Z}_{3}\right)
\rtimes _{\beta }\mathbb{Z}_{2}$ where $\beta :a\in A\mapsto \alpha _{\tau
}(a)$ and $\beta (U_{\eta })=U_{\tau \eta \tau }=U_{\eta }^{2}$. Since $%
A\rtimes _{\alpha _{\eta }}\mathbb{Z}_{3}=A+AU_{\eta }+AU_{\eta }^{2}$, the
relation between $\beta $ and $\alpha _{\eta }$ is given by:%
\begin{equation*}
\beta \left( x_{1}+x_{2}U_{\eta }+x_{3}U_{\eta }^{2}\right) =\alpha _{\tau
}(x_{1})+\alpha _{\tau }(x_{3})U_{\eta }+\alpha _{\tau }(x_{2})U_{\eta }^{2}
\end{equation*}%
for all $x_{1}$,$x_{2}$ and $x_{3}\in A$. We now proceed with a careful
analysis of $\beta $ and $\alpha _{\eta }$ to describe all irreducible
representations of $A\rtimes _{\alpha }\mathfrak{S}_{3}$.

Let $\Pi $ be an irreducible representation of $A\rtimes _{\alpha }\mathfrak{%
S}_{3}$ on some Hilbert space $\mathcal{H}$. Thus $\Pi $ is an irreducible
representation of $\left[ A\rtimes _{\alpha _{\eta }}\mathbb{Z}_{3}\right]
\rtimes _{\beta }\mathbb{Z}_{2}$. We now have two cases: either $\Pi
_{|A\rtimes _{\alpha _{\eta }}\mathbb{Z}_{3}}$ is irreducible or it is
reducible.

\begin{description}
\item[Case 1: $\Pi _{|A\rtimes _{\protect\alpha _{\protect\eta }}\mathbb{Z}%
_{3}}$ is irreducible.] Hence $\Pi $ is minimal for the action $\beta $ of $%
\mathbb{Z}_{2}$. This case splits in two cases.

\begin{description}
\item[Case 1a: $\protect\pi _{A}$ is irreducible] Then $\Pi $ is minimal for
the action $\alpha $ of $\mathfrak{S}_{3}$ by definition.

\item[Case 1b: $\protect\pi _{A}$ is reducible] By Theorem (\ref%
{CyclicConclusion}), there exists an irreducible representation $\pi _{1}$
of $A$ on some Hilbert space $\mathcal{H}_{1}$ such that $\pi _{1}$, $\pi
_{1}\circ \alpha _{\eta }$ and $\pi _{1}\circ \alpha _{\eta }^{2}$ are not
unitarily equivalent, $\mathcal{H}=\mathcal{H}_{1}\oplus \mathcal{H}%
_{1}\oplus \mathcal{H}_{1}$ and:%
\begin{equation*}
\Pi (U_{\eta })=\left[ 
\begin{array}{ccc}
0 & 1 & 0 \\ 
0 & 0 & 1 \\ 
1 & 0 & 0%
\end{array}%
\right] \text{ and\ }\Pi \left( a\right) =\left[ 
\begin{array}{ccc}
\pi _{1}(a) &  &  \\ 
& \pi _{1}\circ \alpha _{\eta }(a) &  \\ 
&  & \pi _{1}\circ \alpha _{\eta ^{2}}(a)%
\end{array}%
\right] \text{.}
\end{equation*}
\end{description}

\item[Case 2: $\Pi _{|A\rtimes _{\protect\beta }\mathbb{Z}_{3}}$ is
reducible.] From Theorem (\ref{CyclicConclusion}), or alternatively \cite%
{Latremoliere06}, there exists an irreducible representation $\pi _{1}$ of $%
A\rtimes _{\alpha _{\eta }}\mathbb{Z}_{3}$ such that for all $z\in A\rtimes
_{\alpha _{\eta }}\mathbb{Z}_{3}$ we have:%
\begin{equation}
\Pi (a)=\left[ 
\begin{array}{cc}
\pi _{1}(z) & 0 \\ 
0 & \pi _{1}\circ \beta (z)%
\end{array}%
\right] \text{ and }\Pi \left( U_{\tau }\right) =\left[ 
\begin{array}{cc}
0 & 1 \\ 
1 & 0%
\end{array}%
\right]  \label{Case2-1}
\end{equation}%
where $\pi _{1}$ and $\pi _{1}\circ \beta $ are not unitarily equivalent.

This case splits again in two cases:

\begin{description}
\item[Case 2a:\ $\protect\pi _{1|A}$ is irreducible] Thus $\pi _{1}$ is a
minimal representation of $A\rtimes _{\alpha _{\eta }}\mathbb{Z}_{3}$. In
particular: 
\begin{equation*}
\pi _{A}(a)=\left[ 
\begin{array}{cc}
\pi _{1}(a) & 0 \\ 
0 & \pi _{1}\circ \alpha _{\tau }(a)%
\end{array}%
\right]
\end{equation*}%
and $\Pi \left( U_{\eta }\right) $ is a block-diagonal unitary in this
decomposition. However, we can not conclude that $\pi _{1|A}$ and $\pi
_{1|A}\circ \alpha _{\tau }$ are equivalent or non-equivalent. Examples (\ref%
{ExPermutation2})\ and (\ref{ExPermutation1})\ illustrate that both
possibilities occur.

\item[Case 2b: $\protect\pi _{1|A}$ is reducible] Then $\Pi _{|A\rtimes
_{\alpha _{\eta }}\mathbb{Z}_{3}}$ is described by Theorem (\ref%
{CyclicConclusion}). Since $3$ is prime, only one possibility occurs: there
exists an irreducible representation $\pi $ of $A$ such that $\pi \circ
\alpha _{\eta }$ and $\pi \circ \alpha _{\eta }^{2}$ are not equivalent and:%
\begin{equation*}
\Pi (a)=\left[ 
\begin{array}{ccc}
\pi (a) & 0 & 0 \\ 
0 & \pi (\alpha _{\eta }(a)) & 0 \\ 
0 & 0 & \pi \left( \alpha _{\eta ^{2}}(a)\right)%
\end{array}%
\right]
\end{equation*}%
and $\Pi \left( U_{\eta }\right) =\left[ 
\begin{array}{ccc}
0 & 1 & 0 \\ 
0 & 0 & 1 \\ 
1 & 0 & 0%
\end{array}%
\right] $. Note that:%
\begin{equation*}
\Pi \left( \beta (U_{\Pi }\right) )=\left[ 
\begin{array}{ccc}
0 & 0 & 1 \\ 
1 & 0 & 0 \\ 
0 & 1 & 0%
\end{array}%
\right] \text{.}
\end{equation*}
Together with (\ref{Case2-1}), we get that $\mathcal{H}$ splits into the
direct sum of six copies of the Hilbert space on which $\pi $ acts and:%
\begin{equation*}
\Pi (a)=\left[ 
\begin{array}{cccccc}
\pi (a) &  &  &  &  &  \\ 
& \pi (\alpha _{\eta }(a)) &  &  &  &  \\ 
&  & \pi (\alpha _{\eta ^{2}}(a)) &  &  &  \\ 
&  &  & \pi \left( \alpha _{\tau }(a)\right) &  &  \\ 
&  &  &  & \pi \left( \alpha _{\eta \tau }(a)\right) &  \\ 
&  &  &  &  & \pi \left( \alpha _{\eta ^{2}\tau (a)}\right)%
\end{array}%
\right]
\end{equation*}%
and%
\begin{equation*}
\Pi (U_{\eta })=\left[ 
\begin{array}{cccccc}
0 & 1 & 0 & 0 & 0 & 0 \\ 
0 & 0 & 1 & 0 & 0 & 0 \\ 
1 & 0 & 0 & 0 & 0 & 0 \\ 
0 & 0 & 0 & 0 & 0 & 1 \\ 
0 & 0 & 0 & 1 & 0 & 0 \\ 
0 & 0 & 0 & 0 & 1 & 0%
\end{array}%
\right]
\end{equation*}%
while%
\begin{equation*}
\Pi \left( U_{\tau }\right) =\left[ 
\begin{array}{cccccc}
0 & 0 & 0 & 1 & 0 & 0 \\ 
0 & 0 & 0 & 0 & 1 & 0 \\ 
0 & 0 & 0 & 0 & 0 & 1 \\ 
1 & 0 & 0 & 0 & 0 & 0 \\ 
0 & 1 & 0 & 0 & 0 & 0 \\ 
0 & 0 & 1 & 0 & 0 & 0%
\end{array}%
\right] \text{.}
\end{equation*}%
Thus $\Pi $ is regular induced by $\pi $, and therefore, as $\Pi $ is
irreducible, $\pi \circ \alpha _{\sigma }$ is not equivalent to $\pi $ for
any $\sigma \in \mathfrak{S}_{3}\backslash \{\limfunc{Id}\}$ by Theorem (\ref%
{RegularIrred}).
\end{description}
\end{description}

This concludes our proof.
\end{proof}

\bigskip We show that all four possibilities above do occur in a nontrivial
manner. We use the generators $\tau $ and $\eta $ as defined in Theorem (\ref%
{Permutation3}). Denote by $e$ the identity of $\left\{ 1,2,3\right\} $.
Notice that $\tau ^{2}=\eta ^{3}=e$ and $\tau \eta \tau =\eta ^{2}$ and $%
\tau \eta ^{2}\tau =\tau $, while:%
\begin{equation*}
\mathfrak{S}_{3}=\left\{ e,\eta ,\eta ^{2},\tau ,\eta \tau ,\eta ^{2}\tau
\right\} \text{.}
\end{equation*}%
In particular, $\left\{ 1,\eta ,\eta ^{2}\right\} $ is a normal subgroup of $%
\mathfrak{S}_{3}$. Now, consider the universal C*-algebra of the free group
on three generators $A=C^{\ast }\left( \mathbb{F}_{3}\right) $ and denote by 
$U_{1},U_{2}$ and $U_{3}$ its three canonical unitary generators. Then we
define the action $\alpha $ of $\mathfrak{S}_{3}$ on $A$ by setting $\alpha
_{\sigma }\left( U_{i}\right) =U_{\sigma (i)}$ for any $\sigma \in \mathfrak{%
S}_{3}$. We now show that this simple example admits in a nontrivial way all
types of representations described in Theorem (\ref{Permutation3}).

\begin{example}
There exists a nontrivial irreducible representation $\pi :C^{\ast }\left( 
\mathbb{F}_{3}\right) \rightarrow M_{2}\left( \mathbb{C}\right) $ such that $%
\pi $ and $\pi \circ \alpha _{\tau }$ are unitarily equivalent, but $\pi $
and $\pi \circ \alpha _{\eta }$ are not. Indeed, set:%
\begin{equation*}
\pi \left( U_{1}\right) =%
\begin{bmatrix}
0 & 1 \\ 
1 & 0%
\end{bmatrix}%
\qquad \pi \left( U_{2}\right) =%
\begin{bmatrix}
0 & -1 \\ 
-1 & 0%
\end{bmatrix}%
\qquad \pi \left( U_{3}\right) =%
\begin{bmatrix}
1 & 0 \\ 
0 & -1%
\end{bmatrix}%
.
\end{equation*}%
We check easily that $\pi $ is an irreducible $\ast $-representation. Since%
\begin{equation*}
\begin{bmatrix}
1 & 0 \\ 
0 & -1%
\end{bmatrix}%
\left[ \pi \circ \alpha _{\tau }\right] 
\begin{bmatrix}
1 & 0 \\ 
0 & -1%
\end{bmatrix}%
=\pi ,
\end{equation*}%
$\pi $ and $\pi \circ \alpha _{\tau }$ are unitarily equivalent. To see that 
$\pi $ and $\pi \circ \alpha _{\eta }$ are not unitarily equivalent, notice
that $\pi \left( U_{1}U_{2}-U_{2}U_{1}\right) =0$ but that:%
\begin{equation*}
\pi \left( U_{2}U_{3}-U_{3}U_{2}\right) =%
\begin{bmatrix}
0 & 2 \\ 
-2 & 0%
\end{bmatrix}%
\text{.}
\end{equation*}
\end{example}

\begin{example}
\label{Minimal}There exists a non trivial irreducible representation $\pi
:C^{\ast }\left( \mathbb{F}_{3}\right) \rightarrow M_{3}\left( \mathbb{C}%
\right) $ such that $\pi $ and $\pi \circ \alpha _{\tau }$ are unitarily
equivalent and $\pi $ and $\pi \circ \alpha _{\eta }$ are also unitarily
equivalent. Let $\lambda =\exp \left( \frac{1}{3}2i\pi \right) $. Define 
\begin{equation*}
\pi \left( U_{1}\right) =\left[ 
\begin{array}{cc}
0 & \lambda \\ 
\lambda ^{2} & 0%
\end{array}%
\right] ,\ \ \ \pi \left( U_{2}\right) =\left[ 
\begin{array}{cc}
0 & \lambda ^{2} \\ 
\lambda & 0%
\end{array}%
\right] \text{ and }\pi \left( U_{3}\right) =\left[ 
\begin{array}{cc}
0 & 1 \\ 
1 & 0%
\end{array}%
\right] \text{.}
\end{equation*}%
Let $V=\left[ 
\begin{array}{cc}
1 & 0 \\ 
0 & \lambda ^{2}%
\end{array}%
\right] $. We check that $V\pi \left( U_{i}\right) V^{\ast }=\pi \left(
U_{\left( i+1\right) \func{mod}3}\right) .$ Then let $W=\pi (U_{3})$. Then $%
W\pi \left( U_{1}\right) W^{\ast }=\pi \left( U_{2}\right) $, $W\pi \left(
U_{2}\right) W^{\ast }=\pi \left( U_{1}\right) $, and $W\pi \left(
U_{3}\right) W^{\ast }=\pi \left( U_{3}\right) $. Thus $\pi $ is a minimal
representation of $C^{\ast }\left( \mathbb{F}_{3}\right) $ for the action $%
\alpha $ of $\mathfrak{S}_{3}$.
\end{example}

\begin{example}
\label{ExPermutation2}There exists an irreducible representation $\pi
:C^{\ast }\left( \mathbb{F}_{3}\right) \rightarrow M_{3}\left( \mathbb{C}%
\right) $ such that $\pi $ and $\pi \circ \alpha _{\eta }$ are unitarily
equivalent, but $\pi $ and $\pi \circ \alpha _{\tau }$ are not: Let $\lambda
=\exp \left( \frac{1}{3}2\pi i\right) $ and define unitaries $T$ and $V$ by 
\begin{equation*}
T=%
\begin{bmatrix}
0 & -\frac{4}{5} & {\Large -}\frac{3}{5} \\ 
\frac{4}{5} & -\frac{9}{25} & \frac{12}{25} \\ 
\frac{3}{5} & \frac{12}{25} & -\frac{16}{25}%
\end{bmatrix}
\text{ \ \ \ and}\ \ \ V=%
\begin{bmatrix}
1 & 0 & 0 \\ 
0 & \lambda & 0 \\ 
0 & 0 & \lambda ^{2}%
\end{bmatrix}%
\end{equation*}

Define 
\begin{equation*}
\pi \left( U_{1}\right) =VTV^{2}\qquad \pi \left( U_{2}\right)
=V^{2}TV\qquad \pi \left( U_{3}\right) =T\text{.}
\end{equation*}

It is clear that $\pi $ and $\pi \circ \alpha _{\eta }$ are unitarily
equivalent. We will show that $\pi $ and $\pi \circ \alpha _{\tau }$ are not
unitarily equivalent. Suppose on the contrary that they are. Then there
exists a unitary $W$ such that $W=W^{\ast }=W^{-1}$ and 
\begin{equation*}
WTW=T,\qquad W\left( VTV^{2}\right) W=V^{2}TV\qquad W\left( V^{2}TV\right)
W=VTV^{2}.
\end{equation*}%
From here we conclude that $VWV$ performs the same transformations, that is 
\begin{eqnarray*}
\left( VWV\right) T\left( VWV\right) ^{\ast } &=&T, \\
\left( VWV\right) \left[ VTV^{2}\right] \left( VWV\right) ^{\ast }
&=&V^{2}TV, \\
\left( VWV\right) \left[ V^{2}TV\right] \left( VWV\right) ^{\ast }
&=&VTV^{2}.
\end{eqnarray*}%
Indeed, 
\begin{eqnarray*}
W\left( VTV^{2}\right) W &=&V^{2}TV\text{ so} \\
V\left[ W\left( VTV^{2}\right) W\right] &=&V\left[ V^{2}TV\right] =TV\text{.}
\end{eqnarray*}
Then we multiply both sides by $V^{2}$ from the right to get 
\begin{equation*}
VWVTV^{2}WV^{2}=T \text{.}
\end{equation*}
Since 
\begin{equation*}
\left( VWV\right) ^{\ast }=V^{\ast }W^{\ast }V^{\ast }=V^{2}WV^{2}\text{,}
\end{equation*}%
we get the first equation. Similarly we get the other two.

Since $\pi $ is irreducible we conclude that there exists a constant $c$
such that 
\begin{equation*}
VWV=cW.
\end{equation*}%
$V$ has a precise form and when we compute $VWV-cW$ we conclude that this
equation has a non zero solution iff $c=1,$ $c=\lambda ,$ or $c=\lambda
^{2}. $ Moreover, the solutions have the form:%
\begin{eqnarray*}
W &=&%
\begin{bmatrix}
x & 0 & 0 \\ 
0 & 0 & y \\ 
0 & z & 0%
\end{bmatrix}%
\text{ if }c=1 \\
W &=&%
\begin{bmatrix}
0 & x & 0 \\ 
y & 0 & 0 \\ 
0 & 0 & z%
\end{bmatrix}%
\text{ if }c=\lambda \\
W &=&%
\begin{bmatrix}
0 & 0 & x \\ 
0 & y & 0 \\ 
z & 0 & 0%
\end{bmatrix}%
\text{ if }c=\lambda ^{2}
\end{eqnarray*}%
for some $x,y,c\in 
%TCIMACRO{\U{2102} }%
%BeginExpansion
\mathbb{C}
%EndExpansion
$.

Now we easily check that $T$ does not commute with any of the three $W$'s.
For example,%
\begin{eqnarray*}
&&%
\begin{bmatrix}
x & 0 & 0 \\ 
0 & 0 & y \\ 
0 & z & 0%
\end{bmatrix}%
\begin{bmatrix}
0 & -\frac{4}{5} & -\frac{3}{5} \\ 
\frac{4}{5} & -\frac{9}{25} & \frac{12}{25} \\ 
\frac{3}{5} & \frac{12}{25} & -\frac{16}{25}%
\end{bmatrix}%
-%
\begin{bmatrix}
0 & -\frac{4}{5} & -\frac{3}{5} \\ 
\frac{4}{5} & -\frac{9}{25} & \frac{12}{25} \\ 
\frac{3}{5} & \frac{12}{25} & -\frac{16}{25}%
\end{bmatrix}%
\begin{bmatrix}
x & 0 & 0 \\ 
0 & 0 & y \\ 
0 & z & 0%
\end{bmatrix}
\\
&=&%
\begin{bmatrix}
0 & \frac{3}{5}z-\frac{4}{5}x & \frac{4}{5}y-\frac{3}{5}x \\ 
\frac{3}{5}y-\frac{4}{5}x & \frac{12}{25}y-\frac{12}{25}z & -\frac{7}{25}y
\\ 
\frac{4}{5}z-\frac{3}{5}x & \frac{7}{25}z & \frac{12}{25}z-\frac{12}{25}y%
\end{bmatrix}%
\text{.}
\end{eqnarray*}%
This of course implies that $x=y=z=0$.
\end{example}

\begin{example}
\label{Torus1}This example acts on $A=C\left( \mathbb{T}^{3}\right) $.
Define for $f\in C(\mathbb{T}^{3})$ and $\left( z_{1},z_{2},z_{3}\right) \in 
\mathbb{T}^{3}$: 
\begin{equation*}
\alpha _{\eta }\left( f\right) \left( z_{1},z_{2},z_{3}\right) =f\left(
z_{2},z_{3},z_{1}\right)
\end{equation*}%
and 
\begin{equation*}
\alpha _{\tau }\left( f\right) \left( z_{1},z_{2},z_{3}\right) =f\left(
z_{2},z_{1},z_{3}\right)
\end{equation*}%
on $C\left( \mathbb{T}^{3}\right) $. We can build a non trivial irreducible
representation $\pi :C\left( \mathbb{T}^{3}\right) \rightarrow 
%TCIMACRO{\U{2102} }%
%BeginExpansion
\mathbb{C}
%EndExpansion
$ such that $\pi $ and $\pi \circ \alpha _{\eta }$ are not unitarily
equivalent and $\pi $ and $\pi \circ \alpha _{\tau }$ are also not unitarily
equivalent. Let $x=\left( x_{1},x_{2},x_{3}\right) \in \mathbb{T}^{3}$ be
such that $x_{1}\neq x_{2},$ $x_{2}\neq x_{3},$ and $x_{3}\neq x_{1}.$Define 
$\pi (f)=f(x)$. Then we obtain an irreducible representation of the required
type as the regular representation induced by $\pi $, using Theorem\ (\ref%
{RegularIrred}).
\end{example}

\bigskip Now, Theorem (\ref{FiniteGroupConclusion}) allowed for the
irreducible subrepresentations of $\Pi _{|A}$ to have multiplicity greater
than one, for irreducible representations $\Pi $ of $A\rtimes _{\alpha }G$.
This situation is however prohibited when $G$ is finite cyclic by Theorem (%
\ref{CyclicConclusion}). We show that finite polycyclic groups such as $%
\mathfrak{S}_{3}$ can provide examples where $\Pi _{|A}$ may not be
multiplicity free, thus showing again that Theorem (\ref%
{FiniteGroupConclusion}) can not be strengthened to the conclusion of
Theorem (\ref{CyclicConclusion}).

\begin{example}
\label{ExPermutation1}We shall use the notations of Theorem (\ref%
{Permutation3}). There exists a unital C*-algebra $A$, an action $\alpha $
of $\mathfrak{S}_{3}$ on $A$ and an irreducible representation $\widetilde{%
\Pi }:A\rtimes _{\alpha }S_{3}\rightarrow B\left( \mathcal{H}\oplus \mathcal{%
H}\right) $ such that for all $x\in A$ we have:%
\begin{equation}
\widetilde{\Pi }\left( x\right) =%
\begin{bmatrix}
\pi \left( x\right) & 0 \\ 
0 & \pi \left( \alpha _{\tau }(x)\right)%
\end{bmatrix}
\label{ExampleMultiplicity2-1}
\end{equation}%
for some irreducible representation $\pi :A\rightarrow B\left( \mathcal{H}%
\right) $ such that $\pi $ and $\pi \circ \alpha _{\tau }$ are equivalent.
Note that $\pi $ is thus minimal for the action of $\alpha _{\eta }$.

Indeed, let us start with any unital C*-algebra $A$ for which there exists
an action $\alpha $ of $\mathfrak{S}_{3}$ and an irreducible representation $%
\Pi :A\rtimes _{\alpha }\mathfrak{S}_{3}\rightarrow B\left( \mathcal{H}%
\right) $ such that $\pi =\Pi _{|A}$ is also irreducible, i.e. $\Pi $ is
minimal. For instance, Example (\ref{Minimal})\ provides such a situation.
Let $V_{\eta }=\Pi \left( U_{\eta }\right) $ and $V_{\tau }=\Pi (U_{\tau })$%
. Then for all $x\in A$%
\begin{eqnarray*}
V_{\eta }\pi \left( x\right) V_{\eta }^{\ast } &=&\pi \left( \alpha _{\eta
}\left( x\right) \right) \text{,} \\
V_{\eta }\pi \left( x\right) V_{\eta }^{\ast } &=&\pi \left( \alpha _{\eta
}\left( x\right) \right) \text{,} \\
V_{\tau }^{2}=1\text{, }V_{\eta }^{3} &=&1\text{ }\text{and }V_{\tau
}V_{\eta }V_{\tau }=V_{\eta }^{2}\text{.}
\end{eqnarray*}%
Let $\omega =\exp \left( \frac{1}{3}2\pi i\right) $. For $x\in A$ define $%
\widetilde{\Pi }\left( x\right) $ by (\ref{ExampleMultiplicity2-1}); let $%
W_{\eta }=\widetilde{\Pi }\left( U_{\eta }\right) $ and $W_{\tau }=%
\widetilde{\Pi }\left( U_{\tau }\right) $ given by:%
\begin{equation*}
W_{\eta }=%
\begin{bmatrix}
\omega V_{\eta } & 0 \\ 
0 & \omega ^{2}V_{\eta }%
\end{bmatrix}%
\text{ \qquad }W_{\tau }=%
\begin{bmatrix}
0 & 1 \\ 
1 & 0%
\end{bmatrix}%
\text{.}
\end{equation*}
We easily check that: 
\begin{equation*}
W_{\eta }\widetilde{\Pi }\left( x\right) W_{\eta }^{\ast }=\widetilde{\Pi }%
\left( \alpha _{\eta }\left( x\right) \right) \text{, } W_{\tau }\widetilde{%
\Pi }\left( x\right) W_{\tau }^{\ast }=\widetilde{\Pi } \left( \alpha _{\tau
}\left( x\right) \right) \text{,}
\end{equation*}
and: 
\begin{equation*}
\left( W_{\tau }\right) ^{3}=1, \left( W_{\tau }\right) ^{2}=1 \text{.}
\end{equation*}
Moreover, 
\begin{eqnarray*}
W_{\tau }W_{\eta }W_{\tau } &=&%
\begin{bmatrix}
0 & 1 \\ 
1 & 0%
\end{bmatrix}%
\begin{bmatrix}
\omega V_{\eta } & 0 \\ 
0 & \omega ^{2}V_{\eta }%
\end{bmatrix}%
\begin{bmatrix}
0 & 1 \\ 
1 & 0%
\end{bmatrix}
\\
&=&%
\begin{bmatrix}
\omega ^{2}\left( V_{\eta }\right) ^{2} & 0 \\ 
0 & \omega \left( V_{\eta }\right) ^{2}%
\end{bmatrix}%
=\left( V_{\eta }\right) ^{2}\text{,}
\end{eqnarray*}%
because $\omega ^{4}=\omega .$

We need to prove that $\widetilde{\Pi }:A\rtimes _{\alpha }S_{3}\rightarrow
B\left( \mathcal{H}\oplus \mathcal{H}\right) $ is irreducible. Let 
\begin{equation*}
T=%
\begin{bmatrix}
a & b \\ 
c & d%
\end{bmatrix}%
\end{equation*}%
be in the commutant of $\widetilde{\Pi }\left( A\rtimes _{\alpha
}S_{3}\right) $. For every $x\in A$: 
\begin{equation*}
T\left[ 
\begin{array}{cc}
\pi (x) & 0 \\ 
0 & \pi (\alpha _{\tau }(x))%
\end{array}%
\right] =\left[ 
\begin{array}{cc}
\pi (x) & 0 \\ 
0 & \pi \left( \alpha _{\tau }(x)\right)%
\end{array}%
\right] T\text{.}
\end{equation*}%
Since $\pi $ is an irreducible representation of $A$ and $\pi \circ \alpha
_{\tau }=V_{\tau }\pi V_{\tau }$ by construction, we conclude by Lemma (\ref%
{Schur})\ that $a$ and $b$ are multiple of the identity, while $c$ and $d$
are multiples of $V_{\tau }$. Since $TW_{\tau }=W_{\tau }T$ we conclude that 
$a=d$ and $b=c$. This means that 
\begin{equation*}
T-aI=%
\begin{bmatrix}
0 & bV_{\tau } \\ 
bV_{\tau } & 0%
\end{bmatrix}%
\end{equation*}%
is in the commutant of the $\widetilde{\pi }\left( A\rtimes _{\alpha
}S_{3}\right) $. However, this element must commute with $W_{\eta }$. This
can only happen if $b=0$. This completes the proof.
\end{example}

\bigskip Thus, using Example (\ref{ExPermutation1}), there exists an
irreducible representation $\widetilde{\Pi }$ of $C^{\ast }(\mathbb{F}%
_{3})\rtimes _{\alpha }\mathfrak{S}_{3}$ such that $\widetilde{\Pi }%
_{|C^{\ast }(\mathbb{F}_{3})}$ is the sum of two equivalent irreducible
representations of $C^{\ast }(\mathbb{F}_{3})$, a situation which is
impossible for crossed-product by finite cyclic groups by Theorem (\ref%
{CyclicConclusion}).

\bigskip In general, repeated applications of Theorem (\ref{CyclicConclusion}%
) can lead to detailed descriptions of irreducible representations of
crossed-products of unital C*-algebra by finite polycyclic groups, based
upon the same method as we used in Theorem (\ref{Permutation3}). Of course,
in these situations Theorem (\ref{FiniteGroupConclusion})\ provides already
a detailed necessary condition on such representations, and much of the
structure can be read from this result.

\providecommand{\bysame}{\leavevmode\hbox to3em{\hrulefill}\thinspace}
\providecommand{\MR}{\relax\ifhmode\unskip\space\fi MR }
\providecommand{\MRhref}[2]{%
  \href{http://www.ams.org/mathscinet-getitem?mr=#1}{#2}
}
\providecommand{\href}[2]{#2}

\end{document}